\documentclass[11pt]{extarticle}

\usepackage{amsfonts}
\usepackage{graphicx} 
\usepackage{color}
\usepackage{subfigure}
\usepackage{amsmath} 
\usepackage{amssymb} 
\usepackage{mathrsfs}
\usepackage{alltt}
\usepackage{cite}
\usepackage{multicol}
\usepackage{multirow}
\usepackage{threeparttable}
\usepackage{float} 

\usepackage{algorithm}
\usepackage{algorithmic}

\usepackage{amsthm}
\usepackage{tabulary}
\usepackage{mathtools}
\usepackage{booktabs}

\usepackage{geometry} 
\newtheorem{Theo}{Theorem}[section]
\newcommand{\bo}{\begin{Theo}}
\newcommand{\eo}{\end{Theo}}
\newtheorem{Lem}{Lemma}
\newcommand{\bl}{\begin{Lem}}
\newcommand{\el}{\end{Lem}}

\newcommand{\E}{\mathbb{E}}            

\def\vec{\operatorname{vec}}

\newtheorem{theorem}{Theorem}[section]

\newtheorem{lemma}[theorem]{Lemma}
\newtheorem{corollary}[theorem]{Corollary}

\newtheorem{assumption}[theorem]{Assumption}
\newtheorem{remark}[theorem]{Remark}

\title{Toward Optimal Statistical Inference in Noisy Linear Quadratic Reinforcement Learning over a Finite Horizon}

\author{Bo Pan$^{a}$, Jianya Lu$^{b}$, Yafei Wang$^{a}$, Hao Li$^{c}$, Bei Jiang$^{a}$, Linglong Kong$^{a}$ \footnote{Corresponding author lkong@ualberta.ca}\\
$^a$Department of Mathematical and Statistical Sciences, University of Alberta.\\
$^b$School of Mathematics, Statistics and Actuarial Science, University of Essex. \\
$^c$School of Economics, LEBPS, Nankai University.}
\date{}

\begin{document}

\maketitle

\begin{abstract}
	
Recent developments in Reinforcement learning have significantly enhanced sequential decision-making in uncertain environments. Despite their strong performance guarantees, most existing work has focused primarily on improving the operational accuracy of learned control policies and the convergence rates of learning algorithms, with comparatively little attention to uncertainty quantification and statistical inference. Yet, these aspects are essential for assessing the reliability and variability of control policies, especially in high-stakes applications such as personalized medicine and optimal liquidation. In this paper, we study statistical inference for the policy gradient (PG) method for noisy Linear Quadratic Reinforcement learning (LQ RL) over a finite time horizon, where linear dynamics with both known and unknown drift parameters are controlled subject to a quadratic cost. In particular, we establish the theoretical foundations for statistical inference in LQ RL, deriving exact asymptotics for both the PG estimators and the corresponding objective loss. Furthermore, we introduce a principled inference framework that leverages online bootstrapping to construct confidence intervals for both the learned optimal policy and the corresponding objective losses. The method updates the PG estimates along with a set of randomly perturbed PG estimates as new observations arrive. We prove that the proposed bootstrapping procedure is distributionally consistent and that the resulting confidence intervals achieve both asymptotic and non-asymptotic validity. Notably, our results imply that the quantiles of the exact distribution can be approximated at a rate of $n^{-1/4}$, where $n$ is the number of samples used during the procedure. The proposed procedure is easy to implement and applicable to both offline and fully online settings. Numerical experiments illustrate the effectiveness of our approach across a range of noisy linear dynamical systems.

\end{abstract}


\section{Introduction}
\label{sec:intro}

Reinforcement learning ({RL}) addresses the problem of learning optimal policies of acting in an uncertain and dynamic environment, which has become an important frontier for several application domains and research areas. Despite {RL}'s promising theoretical and practical performance, limitations persist and restrict its feasibility and scalability in real-world applications \cite{dulac2021challenges, ding2020challenges}. The primary obstacle is insufficient awareness and quantification of uncertainty within the agent or the environment \cite{lockwood2022review}. {RL} seeks policies by interacting with the environment through iterations of \textit{exploration} and \textit{exploitation}, where the former corresponds to learning by interacting with the environment and the latter corresponds to optimizing the objective function given the accumulated information. Herein, randomized actions are generated by random policies to explore the environment and facilitate learning. However, at the same time, they could introduce an implicit bias by assigning a positive probability to potentially non-optimal actions \cite{szpruch2024optimal}. There could also be additional sources of uncertainty due to the lack of extrapolation to unseen states \cite{razin2024implicit, keskin2018incomplete}. These facts raise concerns about robustness and adaptability in cases where agents execute the trained {RL} policy, encountering a new environment. Therefore, it is crucial to develop an effective validation strategy and quantify the uncertainty and variability of {RL} policies. 

The existing literature on {RL} has been restricted to regret, convergence rate, and the stability of the learning algorithm \cite{wagenmaker2022first, littman1996generalized, fallah2021convergence}. Studies of uncertainty quantification and validity testing have only recently attracted attention and have proved essential for this subfield \cite{ramprasad2023online, lockwood2022review, shi2023dynamic, zhu2024uncertainty}. We argue that inferential statistics is one of the key components in quantifying and testing these potentially harmful effects. The recent literature proposes sequential testing procedures to build valid hypothesis tests and the construction of confidence regions of unknown parameters for data-driven decision-making \cite{li2021unifying, fang2018online, shi2023dynamic}. This is mainly because of their critical importance in measuring the statistical uncertainty associated with the estimate in applications compared to a single-point estimate. Meanwhile, implementing a policy without statistical verification of its quality is often risky. For example, in a speed control system for automated driving, a new policy is typically tested on a small fraction of traffic and road conditions before deployment \cite{meadows2022linear}. A confidence interval for the policy can be used to make a more reliable decision because an unstable control policy can potentially reduce the safety of drivers. Likewise, in many real-world applications of operations, policies need to be evaluated on real test tracks before being deployed in real applications, and we focus on not just obtaining the optimal policy but also measuring its associated statistical uncertainty, such as A/B testing conducted by pharmaceutical companies to compare a new product/drug with an old standard. 

Although the development of statistical inference for policy estimation is imperative, studying inference procedures for RL is challenging due to, e.g., the limited theoretical progress, the large memory requirement, and the computational complexity compounded by the learning component. While there has been considerable effort and progress in understanding and working with inference problems in Machine Learning ({ML}), relying solely on these analyses for {ML} to solve {RL} problems is inadequate and inappropriate. One reason is the limited understanding of RL algorithms, particularly regarding their asymptotic behaviour and the dynamics of the underlying system, which makes it challenging to directly apply or extend existing machine learning results. This requires additional analytical techniques. Secondly, there is a discrepancy between the focus of uncertainty (the uncertainty of the data or the model) for the former and the underlying stochasticity for the latter, including, e.g., the uncertainty of the behavior and interactions of the environment in which the agents are trained, and the current limitations of the training process (e.g., the choice of algorithms, the choice of models, the compromise between exploring and exploiting). Therefore, it is important to consider this discrepancy when proposing and analyzing the inferential tools for {RL}. Third, existing inferential tools focus mainly on pivotal statistics \cite{lee2024fast}, Gaussian approximation with various estimates of the unknown asymptotic parameter (e.g., asymptotic covariance matrix \cite{zhu2023online}), batch estimators \cite{chen2020statistical}, or in combination with the multiplier bootstrap approach \cite{chen2020robust}. Extending the classical statistical approach to {RL} could be infeasible as it requires large memory and is often computationally inefficient in sequential data scenarios. Thus, we focus on the research objectives of investigating the exact asymptotics of the {RL} policies and seeking an efficient inference procedure proposed for practical purposes.

\subsection{Main Contribution}

Taking advantage of recent progress in {RL} to find the optimal policy, many previous studies have successfully used the policy gradient  ({PG}) method to solve the optimal control policies, including asymptotic and non-asymptotic properties. In this paper, we specifically focus on the statistical inference for the Polyak-Ruppert averaged iterates of the {PG} estimator for {RL} tasks, considering the Linear Quadratic Reinforcement Learning ({LQ RL}), where a linear function of the current state and the action taken describes the dynamics of the system state, subject to quadratic costs \cite{szpruch2024optimal, umenberger2019robust}. We are the \textit{first} to establish the Central Limit Theorem ({CLT}) results and propose an efficient inference procedure using online bootstrap for the {PG} for \textit{model-based} {LQ RL} over a finite time horizon. In addition, we extend theoretical guarantees to zeroth-order {PG} for \textit{model-free} situations encountering a more complex environment. Using online bootstrapping techniques, we propose a principled inference procedure to construct a confidence interval of the optimal policy. The proposed procedure is easy to implement and can be applied to a broader class of stochastic systems containing {RL} tasks. In particular, the theoretical results are valid under standard assumptions from the {RL} literature, and the proposed approach can be applied to offline and fully online tasks.

In summary, our primary contributions are twofold.

\begin{enumerate}
	\item We obtain the asymptotic normality of Polyak-Ruppert averaging iterates of the policy gradient for the {LQ RL} problem under mild conditions, with applications in both cases of known (model-based) and unknown (model-free) dynamic parameters. We emphasize that, to the best of our knowledge, this study is the \textit{first} to establish a central limit theorem for the Polyak-Ruppert averaging iterates of policy gradient methods in reinforcement learning tasks. More broadly, it also constitutes the \textit{first} CLT result for zeroth-order stochastic gradient methods reported in the literature.
	\item We propose an online bootstrap algorithm to construct a confidence interval for an optimal policy by performing randomly perturbed iteration updates. The proposed bootstrap approach is more applicable for {RL} tasks because it does not require storing iterations in memory and can be performed in parallel. From a theoretical perspective, we establish the asymptotic distributional consistency of the proposed bootstrap algorithm. Furthermore, we conduct a non-asymptotic analysis showing that the quantiles of the exact distribution of the Polyak-Ruppert averaged statistics can be approximated at a rate of $n^{1 / 4}$, providing significant insight into the validity of the proposed bootstrap procedure in a finite-sample manner. The numerical experiments under various settings demonstrate the effectiveness of the proposed method in providing a coverage guarantee for noisy linear dynamic systems.
\end{enumerate}

The remainder of the paper is organized as follows. Section \ref{sec:background} introduces the problem statement for the {LQ RL} task and collects essential notations and lemmas throughout the rest of the paper. Section \ref{sec:policyg} introduces the mathematical framework of the policy gradient for {LQ RL} over a finite time horizon. In Section \ref{sec:bootstrap}, we propose a bootstrapping inference method to construct a confidence interval for the obtained optimal policy. The asymptotic normality of the averaging iterations of the policy gradient for the {LQ RL} task and the distributional consistency of the proposed bootstrap algorithm are established in Section \ref{sec:mainresult}. The results of numerical experiments are presented in Section \ref{sec:ne}, and we conclude with a discussion in Section \ref{sec:discussion}. The proofs of the main results are provided in the appendix.

\subsection{Related Work}

The study of {LQ RL} has been one of the most active areas for both the stochastic control and the {RL} communities \cite{basei2022logarithmic, lovatto2021gradient}, especially for the {RL} with finite state-action and model-based {RL} \cite{carmona2019linear}. It serves as a key benchmark for studying learning policy. First, linear dynamics can be widely used in many applications, such as portfolio allocation \cite{wang2020continuous}, power grids \cite{vlahakis2019distributed}, and resource allocation \cite{ yang2020leveraging}. Second, as a particular case within the broader framework of Reinforcement Learning (RL), this structure has become a critical theoretical benchmark for {RL} tasks, especially when dealing with complex dynamical systems \cite{bradtke1992reinforcement, guo2023reinforcement, wang2021exact, yang2019provably}. This framework is frequently assessed by examining the concept of regret of learning, which quantifies the discrepancy between the value function of the implemented policies and the optimal value function that would be possible if the agent had complete knowledge of the environment. Third, linear dynamics is crucial for providing a reasonable approach to non-linear dynamics: if the dynamics are non-linear and complex to analyze, a local expansion can be used to approximate the original problem \cite{giegrich2024convergence}. 

The Policy Gradient is a popular and principal approach for {RL} tasks for various reasons \cite{fazel2018global, sun2021learning, hambly2021policy}: (1) it is easy to implement; (2) the policy gradient is sufficiently scaled for various types of {RL} variants, directly optimizing the parameter of interest; (3) The property of "gradient dominance" guarantees convergence by running a policy gradient and inherently allows for richly parameterized estimation \cite{sun2021learning}. Moreover, there is growing interest in the development of inference tools for gradient-based approaches in machine learning. For example, \cite{chen2021statistical} developed a stochastic gradient algorithm for online decision-making and proposed plug-in estimators to estimate parameter variance. \cite{li2021statistical} established a functional central limit theorem of the averaged iterates of the local stochastic gradient descent and developed two iterative inference methods, plugging and random scaling, where random scaling constructs an asymptotically pivotal statistic for inference. 

In the realm of statistical inference for {RL}, our paper aligns with a body of literature, including studies such as \cite{ramprasad2023online, hao2019bootstrapping, wang2020residual}, which have introduced bootstrapping techniques for bandit learning and {RL} tasks. Specifically, \cite{ramprasad2023online} formulated an online bootstrapping method for linear stochastic approximation under Markovian noise and applied it to both on-policy and off-policy {RL} algorithms. Similarly, \cite{wang2020residual} advanced a perturbation-based bootstrap exploration method in the bandit setting. Under the same setting, \cite{hao2019bootstrapping} utilized multiplier bootstrapping to estimate the upper confidence bound for exploration. However, no prior study \textit{directly} addressed the statistical inference of the policy gradient estimator in the {LQ RL} framework. Although existing approaches could be extended to {RL}, developing an inference procedure for control problems remains a technical challenge due to the dynamic structure of uncertainty. This challenge necessitates the study of exact and asymptotic expressions for the underlying distribution of the learned policy. In addition, classical statistical inference methods, such as the plug-in procedure and the classical bootstrap for approximating the underlying distributions, are inefficient or, even more concerning, not directly applicable to {RL} problems \cite{hao2019bootstrapping, fang2018online, li2021statistical}. These approaches typically require data to be stored in memory, which may not be feasible due to the sequential learning procedure of the {RL} task. In contrast, the proposed method in this work is computationally efficient and provides a principled method to estimate confidence regions, which is also justified by asymptotic theory.

\section{Problem formulation} \label{sec:background}

In this paper, we investigate the Linear Quadratic Reinforcement Learning (LQ RL) over a finite time horizon of length $T$, where $T$ denotes the total number of observation points. While linear dynamics offer a streamlined description of the underlying process, such control problems remain inherently difficult to solve. However, in LQ RL settings, simpler optimality conditions can be derived using algebraic Riccati equations. Due to their tractability, LQ RL models have become central to optimal control theory and have also attracted growing interest from the reinforcement learning (RL) community. Formally, we consider the {LQ RL} problem as follows \cite{hambly2021policy, yang2019provably, basei2022logarithmic},
\begin{equation} \label{lqr}
	\min _{\left\{u_t\right\}_{t=0}^{T-1}} \mathbb{E}\left[\sum_{t=0}^{T-1}\left(x_t^{\top} Q_t x_t+u_t^{\top} R_t u_t\right)+x_T^{\top} Q_T x_T \right], 
\end{equation}
such that for $t=0,1, \cdots, T-1$,
\begin{equation} \label{dp}
	x_{t+1}=A x_t+B u_t+\omega_t, \quad u_t = K_t x_t, \quad x_0 \sim \mathcal{P},
\end{equation}
where $A \in \mathbb{R}^{d \times d}$ and $B \in \mathbb{R}^{d \times p}$ are referred to as system (transition) matrices, $Q_t \in \mathbb{R}^{d \times d}$ and $R_t \in \mathbb{R}^{p \times p}$ are matrices that parameterize the quadratic costs. $x_t \in \mathbb{R}^d$ is the state of the system with the initial state $x_0$ drawn from a distribution $\mathcal{P}$ with finite second moment. $u_t \in \mathbb{R}^p$ is the action taken at time $t$. 
$\omega$ are independent and identically distributed ({i.i.d.}) random noises with zero mean and finite second moment, that is to say, $\mathbb{E}[\omega]=0$ and $\mathbb{E}\left[\omega \omega^\top\right]$ is positive definite and $\omega$ are independent of $x_0$. Note that the expectation in \eqref{lqr} is taken with respect to ({w.r.t}) both $x_0$ and $w_t$, $t=0,\ldots, T-1$. The goal of linear-quadratic reinforcement learning (LQ RL) is to find a policy $\tau$ that, at each time step $t$, produces a control action $u_t = \tau(h_t)$ based on the full observed history of the system, $h_t = \{x_t, u_{t-1}, x_{t-1}, \ldots, u_1, x_1, u_0\}$. The objective is to maximize system performance while minimizing control effort. The state variables $x_t$ and control inputs $u_t$ are assumed to be transformed such that values of $x_t$ closer to zero indicate better system performance, and values of $u_t$ closer to zero correspond to lower control cost or effort. A cost-minimizing policy in this setting is known as the linear-quadratic regulator, $\tau^{\star}(h_t) = K_t^{\star} x_t$, where the policy $K_t^{\star} \in \mathbb{R}^{d \times p}$ is the optimal solution to a system of equations that depend only on $A$, $B$, $\{Q_t\}_{t = 1}^{T}$, and $\{R_t\}_{t = 1}^{T-1}$. We will review these details in a later section.

\subsection{Preliminaries}

We first collect some notations that will be used throughout the rest of the paper. Then, we present some mathematical tools and lemmas necessary for the following analysis. The commonly used notations are summarized in Table \ref{tab:n_dis} for ease of reference.

We use the notation $\widehat{\cdot}$ to denote observed values or estimates computed given real observations, distinguishing them from the corresponding random variables. We denote by $\boldsymbol{u}:=\left(u_0, \ldots, u_{T-1}\right)$, $\boldsymbol{x}:=\left(x_0, \ldots, x_T\right)$, $\boldsymbol{\omega}:=\left(\omega_0, \ldots, \omega_{T-1}\right)$, $\boldsymbol{Q}:=\left(Q_0, \ldots, Q_T\right)$, and $\boldsymbol{R}:=\left(R_0, \ldots, R_{T-1}\right)$, the profile over the decision period $T$. Given a matrix $Z \in \mathbb{R}^{p}\times \mathbb{R}^d$, denote the vectorization of $Z$ by
\begin{equation*}
	\operatorname{vec}\{Z\}=\left(Z_{11},  Z_{21}, \ldots,  Z_{p 1},  Z_{12}, \ldots,  Z_{p 2}, \ldots,  Z_{1 d}, \ldots,  Z_{p d}\right)^{\top}\in \mathbb{R}^{d p}.
\end{equation*}
For continuously differentiable function $f(\cdot): \mathbb{R}^d \rightarrow \mathbb{R}$ and $v, x\in \mathbb{R}^d$, the directional derivative $\nabla_v f(x)$ is defined by $\nabla_v f(x)\ = \ \lim_{\varepsilon \rightarrow 0} (f(x+\varepsilon v)-f(x)) / \varepsilon$. Let $\nabla f(x)\in \mathbb{R}^d$ and $\nabla^2 f(x)\in \mathbb{R}^{d \times d}$ be the gradient and the Hessian matrix of $f$, respectively. We also denote $\otimes$ the Kronecker product. To learn the optimal policy for the {LQ RL} problem as in \eqref{lqr}, we make the following mild assumption on the model parameters to ensure that problem \eqref{lqr} is well defined.

\begin{assumption}[Cost Parameter] \label{assu1}
	Assume $Q_t \in \mathbb{R}^{d \times d}$, for $t=0,1, \ldots, T$, and $R_t \in \mathbb{R}^{p \times p}$, for $t=0,1, \ldots, T-1$, are positive definite matrices.
\end{assumption}

Under Assumption \ref{assu1}, we can properly define a sequence of matrices $\{P^{\star}_t\}_{t=0}^T$ as the solution to the discrete algebraic Riccati equation \cite{bertsekas2012dynamic},
\begin{equation} \label{eq:P1}
	P_t^{\star}=Q_t+A^{\top} P_{t+1}^{\star} A-A^{\top} P_{t+1}^{\star} B\left(B^{\top} P_{t+1}^{\star} B+R_t\right)^{-1} B^{\top} P_{t+1}^{\star} A,
\end{equation}
with the terminal condition $P^{\star}_T=Q_T$. The matrices $\{P^{\star}_t\}_{t=0}^{T-1}$ can then be found by solving equation \eqref{eq:P1} backward iteratively in time, and the solutions $\{P^{\star}_t\}_{t=0}^{T-1}$ are positive definite. In particular, assume that the optimal control for the problem \eqref{lqr} is linearly related to the state of the system. The optimal policy has an explicit formula, as shown in Lemma \eqref{lem::os}.

\begin{lemma}[Optimal Solution of {LQ RL} \cite{bertsekas2012dynamic}] \label{lem::os}
	Under Assumption \ref{assu1}, the optimal control sequence $\left\{u_t\right\}_{t=0}^{T-1}$ is given by
	\begin{equation*} \label{trueK}
		u_t = -K_t^{\star} x_t, \quad \text { where } K_t^{\star} =\left(B^{\top} P_{t+1}^{\star} B+R_t\right)^{-1} B^{\top} P_{t+1}^{\star} A.
	\end{equation*}
\end{lemma}

Lemma \ref{lem::os} indicates that we only need to focus on the linear admissible policies for the {LQ RL} problem in the form of $u_t=-K_t x_t$, $t=0,1, \cdots, T-1$, which can be fully characterized by $\boldsymbol{K}:=(K_1,\ldots, K_{T-1})$ with $K_t =\left(B^{\top} P_{t+1} B+R_t\right)^{-1} B^{\top} P_{t+1} A$, with 
\begin{equation*}
	P_t=Q_t+A^{\top} P_{t+1} A-A^{\top} P_{t+1} B\left(B^{\top} P_{t+1} B+R_t\right)^{-1} B^{\top} P_{t+1} A.
\end{equation*}

When parameters are unknown in model-based {LQ RL}, the classical approach is based on the \textit{certainty equivalence principle} \cite{rathnam2023unintended, aastrom1995adaptive, hambly2021policy}. This method involves using historical data collected under a set policy and applying classical statistical techniques, e.g., least squares estimation \cite{dean2020sample}, to estimate the parameters. Once estimated, these parameters are treated as if they were known with certainty, and then a policy is designed that treats the estimated parameters as the ground truth \cite{mania2019certainty}. For the identification of the system in {LQ RL}, the main ingredient requested is the trajectory of the observations by introducing random noise, creating $m$ replications for each given an initial state and a fixed policy. 
This process allows for the estimation of the system's dynamics under various conditions. Denote the obtained data set as $\left\{\left(\widehat{x}_{t, i}, \widehat{u}_{t, i}\right): 1 \leq i \leq m, 0 \leq t \leq T\right\}$, where $t$ indexes time and $i$ indexes independent replication. Therefore, we can estimate the system dynamics by
\begin{equation*}
	\left(\widehat{A}, \widehat{B}\right) = \arg \min _{(A, B)} \sum_{i=1}^m \sum_{t=0}^{T-1} \frac{1}{2}\left\|A \widehat{x}_{t, i}+B \widehat{u}_{t, i}-\widehat{x}_{t+1, i}\right\|_2^2.
\end{equation*}
Although the ``plug-in" method is an option to construct statistical inference with an estimated covariance structure based on the Gaussian approximation, we emphasize that the bias introduced (approximation error) must be quantified to justify the effectiveness of the plug-in method for statistical inference of the model-based {LQ RL}. Additionally, a fundamental problem of uncertainty quantification is the source of randomness. In {RL}, the source of randomness could come from the problem itself, noise, and randomness from the algorithm. Remarkably, considering the dynamic structure, these random variables have a joint impact on the future states. The approach based on the certainty equivalence principle only considers the uncertainty from the estimation of parameters and ignores the environment noise and stochastics from the algorithm. Following this strategy is known to be inconsistent \cite{becker1985adaptive, lai1982iterated}. Moreover, the covariance matrix contains the inverse of the matrix, which may not exist in practice in a high-dimensional setting. Therefore, it is important to consider the uncertainty of all components of the learning process and their properties, which could differ for different algorithms. 

\begin{assumption} [Stability] \label{assu2}
	Assume that the linear dynamics is stabilizable, i.e. there exists $K$ such that the spectral radius (maximum absolute eigenvalue) of $A+B K$ is strictly less than $1$.
\end{assumption}

\begin{lemma}[\cite{hambly2021policy}] \label{lem::matpos}
	Assume that $x_0$ and $\{\omega_t\}_{t=1}^{T-1}$ have a finite second moment and $\{\omega_t\}_{t = 1}^{T-1}$ are independent of $x_0$. We have that $\mathbb{E}\left[\omega_t \omega_t^\top\right]$ is positive definite for $t=1,\ldots,T-1$ under any policy $\boldsymbol{K}$.
\end{lemma}

\section{Policy gradient descent for {LQ RL}} \label{sec:policyg}

In this section, we introduce the model-based and model-free {PG} for the {LQ RL}. Note that the matrices $\{K_t\}_{t=1}^{T-1}$ can be obtained from the model parameters and $\{P_t\}_{t=1}^{T-1}$ computed via the Riccati equation \eqref{eq:P1}. However, it should be emphasized that this process of obtaining $\{K_t\}_{t=1}^{T-1}$ and $\{P_t\}_{t=1}^{T-1}$ requires computing the inverse of certain matrices. This process may lead to high computational costs and accumulation of computation errors, particularly in high-dimensional scenarios, and can only be used in model-based cases. This paper focuses on using the policy gradient method to solve the {LQ RL} problem \eqref{lqr}.

\subsection{Exact policy gradient descent for model-based {LQ RL}}

We first assume that all parameters of the model $\{Q_t\}_{t = 1}^{T}$, $\{R_t\}_{t = 0}^{T-1}$, $A$, $B$ are known, according to the certainty equivalence principle. When the true underlying dynamics are linear, the model-based approach could outperform others, fully utilizing the linear structure. In addition, the analysis of {PG} with known parameters could provide a learning pathway for the model-free {LQ RL}, bridging zeroth-order {PG}, and inference for the model-free {RL}. 

Since an admissible policy can be fully characterized by $\boldsymbol{K}$, the cost function of a policy can be correspondingly defined as
\begin{equation} \label{e:cost}
	C(\boldsymbol{K})=\mathbb{E}\left[\sum_{t=0}^{T-1}\left(x_t^{\top} Q_t x_t+u_t^{\top} R_t u_t\right)+x_T^{\top} Q_T x_T\right],
\end{equation}
where $\left\{x_t\right\}_{t=1}^T$ and $\left\{u_t\right\}_{t=0}^{T-1}$ are the states and controls induced by $\boldsymbol{K}$, starting with $x_0$. Note that $\boldsymbol{K}^{\star}$ is the optimal policy for the problem \eqref{lqr}, in that
\begin{equation}\label{eq::os}
	\boldsymbol{K}^{\star} = \arg \min_{\boldsymbol{K}}  C(\boldsymbol{K}),
\end{equation}
subject to the dynamics \eqref{dp}. Based on the policy gradient method, the updating rule for the cost function \eqref{e:cost} has the form
\begin{equation} \label{update}
	{K}_t^{\ell+1} = {K}_t^\ell-\eta_\ell \nabla_{t} C\left(\boldsymbol{K}^\ell\right), \quad \forall t=0,1, \ldots, T-1,
\end{equation}
where $\ell$ is the number of iterations with maximum number iteration of $n$, $\nabla_t C(\boldsymbol{K})=\frac{\partial C(\boldsymbol{K})}{\partial K_t}$ is the gradient of $C(\boldsymbol{K})$ {w.r.t} $K_t$, and $\eta_\ell$ is the step size. Lemma \ref{lem:gradient} provides the explicit formulation for the vectorized gradient of the cost function $\nabla_{t} C(\boldsymbol{K})$.

\begin{lemma}[Exact Policy Gradient Descent]\label{lem:gradient}
	For any $t=0,1, \ldots, T-1,$ the policy gradient has the representation of
	\begin{equation*}
		\nabla_{t}C(\boldsymbol{K})=2\operatorname{vec}\left\{{E}_t(K_t)\Sigma_t\right\},
	\end{equation*}
	where ${E}_t(K_t)=\left(R_t+B^{\top} P_{t+1} B\right) K_t-B^{\top} P_{t+1} A$, and $\Sigma_t=\mathbb{E}\left[x_tx_t^\top\right]$. And, $P_t^{\boldsymbol{K}}$ can be defined as the solution to
	\begin{equation} \label{pgP}
		P_t^{\boldsymbol{K}}=Q_t+K_t^{\top} R_t K_t+\left(A-B K_t\right)^{\top} P_{t+1}^{\boldsymbol{K}}\left(A-B K_t\right), \quad t=0,1, \cdots, T-1,
	\end{equation}
	with terminal condition $P_T^{\boldsymbol{K}}=Q_T$. Note that \ref{pgP} is equivalent to the Riccati equation \ref{eq:P1} with optimal $K_t=K_t^{\star}$ given by Lemma \ref{lem::os}. 
\end{lemma}


Lemma \ref{lem:gradient} indicates that the optimal policy for problem \eqref{eq::os} can be obtained by iteratively updating $K^\ell_t$ as follows
\begin{equation}\label{e:sgd}
	\operatorname{vec}\left\{K^{\ell+1}_t\right\}=\operatorname{vec}\left\{K^{\ell}_t\right\}-2 \eta_\ell \operatorname{vec}\left\{{E}_t\left(K_t^{\ell}\right)x_t^{\ell}\left(x_t^\ell\right)^\top\right\}, 
\end{equation}
where $x_t^\ell$ indicates the sample of $x_t$ at $\ell$-th iteration. Given that $x_0$ is sampled from the distribution $\mathcal{P}$, the distribution of the state variable $x_t$, for $t=1,\ldots, T-1$, denoted as $\mathcal{P}_t$, can be derived from the dynamic system \eqref{dp}. Comparing \ref{pgP} and \ref{eq:P1}, we can see that \ref{pgP} does not involve the inverse of matrices, thus avoiding the approximation error and numerical difficulty. Assume $x_{t,i}$, for $i=1, \ldots, m$, are {i.i.d.} random samples from distribution $\mathcal{P}_t$ 
with corresponding actions $u_{t, i}$, for $i = 1,\ldots, m$, determined accordingly. Mirroring the cost function \eqref{e:cost}, we have the empirical risk function given by
\begin{equation*}
	\frac{1}{m}\sum_{i=1}^{m} C(\boldsymbol{K})=\frac{1}{m}\sum_{i=1}^{m}\left\{\sum_{t=0}^{T-1}\left(x_{t,i}^\top Q_t x_{t,i}+u_{t,i}^\top R_t u_{t,i}\right)+x_{T,i}^\top Q_T x_{T,i}\right\}.
\end{equation*}
Algorithm \ref{alg:PG} outlines the iterative procedure used in the policy gradient method, which is tailored for solving the {LQ RL} problem \eqref{e:cost}.

\begin{algorithm}[ht!]
	\caption{Exact policy gradient descent for {LQ RL}}
	\begin{algorithmic}[1]\label{alg:PG}
		\STATE \textbf{Input:} Finite time horizon $T$, maximum number of iterations $n$, initial status $\{x_{0,i}\}_{i = 1}^{m}$, initial policy $K^0$, and stepsize $\left\{\eta_\ell\right\}_{\ell=1}^n$. 
		\FOR{ $t$ in $0, \cdots, T - 1$} 
		\FOR{ $\ell$ in $0, \cdots, n$}
		\STATE Sample $i$ from $\{1, \cdots, m\}$ and compute gradient of $C(\boldsymbol{K})$ as Lemma \ref{lem:gradient}.
		\STATE Update $K_t^{\ell}$: $\operatorname{vec}\{K^{\ell+1}_t\}=\operatorname{vec}\{K^{\ell}_t\}-\eta_\ell \nabla_{t}C(\boldsymbol{K})$, and $P_t^{\boldsymbol{K}}$ as \eqref{pgP}.
		\ENDFOR
		\ENDFOR
		\STATE \textbf{Output:} $\boldsymbol{K} = \{K_t^{n}\}_{t=1}^{T-1}$.
	\end{algorithmic} 
\end{algorithm}

However, similar to other model-based approaches that critically rely on the model assumption, it also suffers from model misspecification. Compared to the model-based approach, on the other hand, the execution of the zeroth-order {PG} for a model-free setting does not rely on the assumptions of the dynamic structure, where the optimal policy was learned directly by interacting with the environment, without inferring the model parameters \cite{pong2018temporal, huang2020model, hambly2021policy}. Thus, the advantage of the zeroth-order {PG} is that it is more robust against model misspecification than the exact {PG}.

\subsection{Zeroth-order policy gradient descent for model-free {LQ RL}}

In this section, we introduce the zeroth-order {PG} for the model-free {LQ RL}. In this model-free setting, the decision maker has no prior information on the model parameters $A$, $B$, $\{Q_t\}_{t = 0}^{T}$ and $\{R_t\}_{t = 0}^{T-1}$ (i.e., assumed to be unknown) but has only access to a sample of the risk function $\widehat{c}(\cdot)$ given arbitrary state $x$ and input $K$. Using the technique of zeroth-order optimization, the gradient of $K$ can be approximated. The zeroth-order {PG} has been shown to converge to the optimal solution in the literature, with high sample efficiency and polynomial computational complexity (see, e.g., \cite{ji2019improved, wang2018stochastic, nikolakakis2022black}). The basic idea of zeroth-order {PG} is to approximate the gradient using the function values as follows

\begin{lemma} [Unbiased Gradient Estimator]\label{lem:gradientfree}
	Under Assumptions \ref{assu1} and \ref{assu2} with the uniform distribution $\mathcal{U}(r)$ with radius $r$, we have
	\begin{eqnarray*}
		\E\left[\vec\left\{\big(c\left(\boldsymbol{K}+\boldsymbol U_t\right)-c\left(\boldsymbol{K}\right)\big) \boldsymbol U_t\right\}\right]= \frac{r^2}{3}\nabla_{t}C(\boldsymbol{K}),
	\end{eqnarray*}	
	where the expectation is taken w.r.t $\mathcal{U}(r)$. Thus, the gradient of the function $C(\boldsymbol{K})$ with input $\boldsymbol{K}$ can be estimated by averaging the observations of $ \frac{3}{r^2}\{{c}(\boldsymbol{K} +\boldsymbol U_t) - {c}(\boldsymbol{K})\} \boldsymbol U_t$.
\end{lemma}


\begin{remark}
	We emphasize that the idea of approximating the gradient by function values is motivated by a first-order Taylor expansion, such that 
	\begin{eqnarray*}
		\E[\big(c(\boldsymbol{K}+\boldsymbol U_t)-c(\boldsymbol{K})\big) \boldsymbol U_t]= \E[\big(c(\boldsymbol{K}+\boldsymbol U_t)\big) \boldsymbol U_t] \approx \mathbb{E}[(\nabla C(\boldsymbol{K}) \boldsymbol U_t) \boldsymbol U_t]= \nabla C(\boldsymbol{K}) L,
	\end{eqnarray*}
	for some constant $L$. In estimating the gradient term as in Lemma \ref{lem:gradientfree}, we only query the access to the observation of $c(\cdot)$ at the input $\boldsymbol{K}$, without information of gradients and higher-order derivation. Our approach is similar to the zeroth-order optimization proposed by \cite{fazel2018global}, where $\mathbb{E}_{\boldsymbol{U} \sim \mathcal{N}\left(0, \sigma^2 I\right)}[c(\boldsymbol{K}+\boldsymbol{U})\boldsymbol{U}]$ is used for estimation. However, when Gaussian smoothing is applied, $c(\cdot)$ could not be bounded for each $\boldsymbol{K}$. This could make the estimator not well-defined. To avoid this, we adopt the approach by generating $\boldsymbol{U}$ from a uniform distribution for each element of $\boldsymbol{K}$. A similar idea can also be found in \cite{hambly2021policy}, where the author proposes to sample $\boldsymbol{U}$ uniformly on the surface of a sphere with a radius of $r$ in $\mathbb{R}^{d \times p}$.
\end{remark}

We denote the policy estimated by zeroth-order {PG} as $\tilde{K}_t$ to distinguish from exact {PG} iterates. With the gradient approximation, we can update the policy in each iteration as follows, 
\begin{equation*}
	\tilde{K}_t^{\ell+1}=K_t^\ell-\eta_\ell \widehat{\nabla}_{t} C\left(\boldsymbol{\tilde{K}}^\ell\right),
\end{equation*}
for $t = 0, \cdots, T-1$, where $\widehat{\nabla}_t C\left(\boldsymbol{\tilde{K}}^\ell\right)$ is the estimate of $\nabla_t C\left(\boldsymbol{\tilde{K}}^\ell\right)$, and the corresponding empirical form can be formulated as follows
\begin{eqnarray} \label{e:sgdfree1}
	\vec\left\{ K_t^{\ell+1}\right\}=\vec\left\{ K_t^{\ell}\right\}-\frac{3\eta_\ell}{r^2m}\sum_{i=1}^{m}\vec\left\{\left(\widehat{c}\left( \tilde{K}_t^\ell+ U_t^i \right)-\widehat{c}\left(\tilde{K}_t^\ell\right)\right)  U_t^i\right\},
\end{eqnarray}
where $U_t^i$ denotes the random matrix at stage $t$, whose elements are independently and uniformly distributed according to $\mathcal{U}(r)$. In each iteration, we randomly generate $m$ random matrices.

\begin{algorithm}[ht!]
	\caption{Zeroth-order policy gradient descent for {LQ RL}}
	\begin{algorithmic}[1]\label{alg:zeroPG}
		\STATE \textbf{Input:} $\boldsymbol{K}$, number of trajectories $m$, parameter $r$.
		\FOR{$i$ in $1, \cdots, m$ }
		\FOR{$t$ in $0, \cdots, T - 1$ }
		\STATE Sample the (sub)-policy at time $t$.
		\STATE Compute ${K}_t^i=K_t+U_t^i$ where $U_t^i$ is drawn independently at random over matrices from the uniform distribution, such that $U_i \sim \mathcal{U}(r)$.
		\STATE Return $\widehat{c}^i(\cdot)$ as the observed single trajectory cost with policy
		\begin{equation*}
			\left(\boldsymbol{K}_{-t}, {K}_t^i\right):=\left(K_0, \cdots, K_{t-1}, {K}_t^i, K_t, \cdots, K_{T-1}\right),
		\end{equation*}
		starting from $x_0^i \sim \mathcal{P}$.
		\ENDFOR
		\ENDFOR
		\STATE Return the gradient approximation of $\nabla_t C(\boldsymbol{K})$ for each $t$
		\begin{equation*}
			\widehat{\nabla}_t C(\boldsymbol{K}) = \frac{3}{r^2m}\sum_{i=1}^{m}\vec\left\{\left(\widehat{c}\left(K^\ell_t+ U_t^i\right)-\widehat{c}\left(K^\ell_t\right)\right) U_t^i\right\} .
		\end{equation*}
	\end{algorithmic} 
\end{algorithm}  

\section{Statistical inference for policy gradient method for {LQ RL}} \label{sec:bootstrap}

As suggested by \cite{ruppert1988efficient, polyak1992acceleration}, we consider the Polyak-Ruppert averaging iterate, 
\begin{equation*}
	\bar{{K}}_t^{n} = \frac{1}{n} \sum_{\ell = 1}^{n} {K}_t^\ell, \quad t=0,1, \cdots, T-1, 
\end{equation*}
where $K_t^\ell$ is defined as in \eqref{e:sgd}. To perform statistical inference for the averaging policy gradient estimator $\bar{{K}}_t^{n}$ at any stage $t$, we first need to establish the asymptotic distribution of $\bar{{K}}_t^{n}$ to show that the study of the inference procedure is appropriate for the {LQ RL} problem. This theoretical property is shown in Section \ref{sec:mainresult}. Unfortunately, it is still challenging to construct statistical inference with an averaging estimator $\bar{{K}}_t^{n}$ because the asymptotic covariance crucial for the inference procedure depends on the unknown value $K^\star_t$. 

To address this difficulty, in this paper, we propose a bootstrapping procedure that updates the policy gradient estimator and a large number of randomly perturbed policy gradient estimators for each iteration step. Specifically, let $W$ be a positive {i.i.d.} random variable with a mean of one and a variance of one. Using non-negative random weights can guarantee a convex-weighted objective function. In parallel to the main update \eqref{update}, the randomly perturbed policy gradient is performed for $\ell = 1, \cdots, n$,     
\begin{equation*}
	K_t^{\ell+1, \prime}=K_t^{\ell, \prime}-\eta_{\ell} W_{\ell} \nabla_{t} C\left(\boldsymbol{K}^{\ell, \prime}\right), \quad t=0,1, \cdots, T-1.
\end{equation*}
Denote $\bar{K}_t^{n, \prime} = \frac{1}{n} \sum_{\ell = 1}^{n} K_t^{\ell, \prime}$, $t=0,1, \cdots, T-1$. We will show that 
\begin{equation*}
	\sqrt{n}\left(\operatorname{vec}\left\{\bar{{K}}_t^{n}\right\} -\operatorname{vec}\left\{K^{\star}_t\right\}\right) \text{ and } \sqrt{n} \left(\operatorname{vec}\left\{\bar{K}_t^{n, \prime}\right\} - \operatorname{vec}\left\{\bar{{K}}_t^n\right\}\right)
\end{equation*}
converge in distribution to the same limiting distribution in Section \ref{sec:mainresult}. It validates the proposed bootstrapping policy gradient procedure's effectiveness, enabling us to conduct inference on the former distribution by using the latter as a proxy. In practice, we generate a large number, say $\mathcal{B}$, of random samples $\{W_{\ell, (b)}\}_{b=1}^\mathcal{B}$ of $W$. For each $b = 1, \cdots, \mathcal{B}$, at time step $\ell$, we update the perturbed policy gradient iteration $K_{t, (b)}^{\ell, \prime}$ as follows,
\begin{equation*}
	K_{t, (b)}^{\ell + 1, \prime} = K_{t, (b)}^{\ell, \prime} - \eta_\ell W_{\ell, (b)} \nabla_{t} C\left(\boldsymbol{K}^{\ell, \prime}\right), \quad \bar{K}_{t, (b)}^{n,\prime} = \frac{1}{n} \sum_{\ell = 1}^n K_{t, (b)}^{\ell, \prime},
\end{equation*}
where $W_{\ell, (b)}$ are {i.i.d.} random variables with mean one and variance one. Parallel to iteration \eqref{e:sgd}, the perturbed bootstrap policy gradient iteration for the model-based and model-free cases are given as
\begin{equation*} \label{e:bsgd}
	\begin{split}
		\operatorname{vec}\left\{K^{\ell+1, \prime}_{t, (b)}\right\} &=\operatorname{vec}\left\{K^{\ell, \prime}_{t, (b)}\right\}-2\eta_\ell W_{\ell, (b)}\operatorname{vec}\left\{E_t\left(K_{t, (b)}^{\ell,\prime}\right)\left[x_t^\ell\left(x_t^\ell\right)^\top\right]\right\},  \\
		\vec\left\{\tilde{K}_{t, (b)}^{\ell+1, \prime}\right\} &=\vec\left\{\tilde{K}_{t, (b)}^{\ell, \prime}\right\} 
		- \frac{3\eta_\ell}{r^2m} W_{\ell, (b)}\sum_{i=1}^{m}\vec\left\{\left(\widehat{c}\left(\tilde{K}^{\ell, \prime}_{t, (b)} + U_{t, (b)}^{i}\right) - \widehat{c}\left(\tilde{K}_{t, (b)}^{\ell, \prime}\right)\right)U_{t, (b)}^{i}\right\}.
	\end{split}   
\end{equation*}
The sampling distribution of $\sqrt{n}\left(\operatorname{vec}\left\{\bar{{K}}_t^{n}\right\} - \operatorname{vec}\left\{K^{\star}_t\right\}\right)$ can then be approximated by the empirical distribution of $\left\{\operatorname{vec}\left\{\bar{K}_{t, (b)}^{n, \prime}\right\}-\operatorname{vec}\left\{\bar{K}_t^{n}\right\}: b=1,\ldots,\mathcal{B}\right\}$. Specifically, the covariance matrix of $\operatorname{vec}\left\{\bar{K}_t^{n}\right\}$ can be estimated by the sample covariance matrix constructed from $$\left\{\operatorname{vec}\left\{\bar{K}_{t, (b)}^{n, \prime}\right\}: b=1,\ldots,\mathcal{B}\right\}.$$ 
It leads to the construction of a confidence interval for $\operatorname{vec}\{K^{\star}_t\}$. We emphasize that the proposed bootstrapping inference procedure retains the simplicity of the policy gradient descent. In addition, all $\mathcal{B}$ trajectories of perturbed iterations depend on a single trajectory of iteration. Therefore, iterates can be updated in a parallel and fully online manner. Algorithm \ref{alg:bpg} summarizes the iteration step for the bootstrapping policy gradient descent method for {LQ RL} to find a confidence interval of $\operatorname{vec}\{K^{\star}_t\}$.   

\begin{algorithm}[ht!]
	\caption{Bootstrapping policy gradient descent for {LQ RL}}
	\begin{algorithmic}[1] \label{alg:bpg}
		\STATE \textbf{Input:} Finite time horizon $T$, maximum number of iterations $n$, initial status $\{x_{0,i}\}_{i = 1}^{m}$, initial policy $\boldsymbol{K}^0$, stepsize $\left\{\eta\right\}_{\ell=1}^n$, number of bootstrap samples $B$.
		\FOR{ $t$ in $0, \cdots, T - 1$} 
		\FOR{ $\ell$ in $0, \cdots, n$}
		\STATE Sample $i$ from $\{1, \cdots, m\}$.
		\STATE Compute gradient of $C(\boldsymbol{K})$ as in Algorithm \ref{alg:PG} or \ref{alg:zeroPG}.
		\STATE Update $K^{\ell+1}_t=\operatorname{vec}\left\{K^{\ell}_t\right\}-\eta_\ell \nabla_{t}C(\boldsymbol{K}^\ell)$.
		\STATE Update $\bar{{K}}_t^{\ell+1} = \frac{1}{\ell+1}\left(\ell\bar{K}_t^{\ell} + K_t^{\ell+1}\right)$.
		\FOR{ $b = 1, \cdots, B$}
		\STATE Update $K_{t, (b)}^{\ell+1, \prime}=K_{t, (b)}^{\ell, \prime}-\eta_\ell W_{\ell, (b)}\nabla_{t} C\left(\boldsymbol{K}^{\ell, \prime}_{(b)}\right)$.
		\STATE Update $\bar{K}_{t, (b)}^{\ell+1, \prime} = \frac{1}{\ell+1}\left(\ell\bar{K}_{t, (b)}^{\ell, \prime} + K_{t, (b)}^{\ell+1, \prime}\right)$.
		\ENDFOR
		
		\ENDFOR
		\ENDFOR
		\STATE Obtain Bootstrap estimate $\left\{\bar{K}_{t, \prime}^{n, \prime}\right\}_{b = 1}^{B}$ and construct confidence interval from the empirical distribution of $\left\{\bar{K}_{t, (b)}^{n, \prime}\right\}_{b = 1}^{B}$, denoted as $[l_t^n, u_t^n]$, $\forall t$.
		\STATE \textbf{Output:} Confidence interval $[l_t^n, u_t^n]$, $\forall t$.
	\end{algorithmic}
\end{algorithm}

\begin{remark}
	Two approaches can be used when we construct the confidence interval based on the empirical distribution of the output of the bootstrap estimate $\left\{\bar{K}_{t, (b)}^{n, \prime}\right\}_{b = 1}^{\mathcal{B}}$. The first approach employs the quantiles of these bootstrap estimates directly. Specifically, define $q_{\delta}$ as the $\delta$-th quantile of the empirical distribution of the bootstrap estimates $\left\{\bar{K}_{t, (b)}^{n, \prime} - \bar{K}_t^n\right\}_{b = 1}^{\mathcal{B}}$, and then a confidence interval can be calculated as $\left[\bar{K}_t^n + q_{\alpha / 2}, \bar{K}_t^n + q_{1 - \alpha/2}\right]$. The second approach relies on the standard error of the covariance matrix. Here, using $z_\alpha$, which denotes the $\alpha$-th quantile of the standard normal distribution, a confidence interval is constructed as $\left[\bar{K}_t^n - z_{1 - \alpha/2}\widehat{\sigma}_t, \bar{K}_t^n + z_{1 - \alpha/2}\widehat{\sigma}_t\right]$, where $\widehat{\sigma}_t = \sqrt{\operatorname{diag}{\widehat{\Sigma}_t}}$, and $\widehat{\Sigma}_t$ is the sample covariance matrix derived from bootstrap estimates.
\end{remark}

\section{Theoretical properties} \label{sec:mainresult}

In this section, we study the theoretical foundation of statistical inference of the policy gradient descent method applied to the {LQ RL} problem. We perform a rigorous analysis to demonstrate that the confidence interval constructed by the proposed bootstrapping algorithm is asymptotically valid. Furthermore, we emphasize that these theoretical results can also be extended to variations of the {LQ RL} model, thus broadening the scope and applicability of our inferential framework. We briefly state the required assumptions and proof sketch in the following.

\begin{assumption}\label{assu3}
	Assume that $x_0\sim\mathcal{P}$ is such that $\E\left[x_0x_0^\top\right]$ is positive definite and $\omega$ are {i.i.d.} random variables with finite second moment bounded and independent of $x_0$ such that $\E[\omega]=0$ and $\E\left[\omega\omega^\top\right]$ are positive definite. 
\end{assumption}

\begin{assumption}\label{assu4}
	The step sizes are chosen as $\eta_\ell = C\eta_\ell^{-\alpha}$ with some $C>0$ and $\alpha\in[1/2,1)$.
\end{assumption}

\begin{assumption}\label{assu5}
	$W$ is a sequence of {i.i.d.} positive random variables with mean and variance equal to one. 
\end{assumption}

The assumptions (Assumption \ref{assu1} -Assumption \ref{assu5}) used in our analysis adhere closely to those commonly accepted within the {LQ RL} literature. The assumptions made are imposed to ensure the stability of the algorithm and dynamics, to ensure convergence with a rate of $\sqrt{n}$ of the sequence of iterations, and to control the behaviour of random noise, which typically holds when the maps and the objective function are bounded. In the literature, polynomially decaying step sizes are also a standard for iterative averaging. We emphasize that, aside from assuming finite second moments of the noise, no distributional assumptions are made, which enhances the model's applicability to real-world scenarios despite the increased complexity this entails.

\subsection{Asymptotic theory for {PG}}

We establish a central limit theorem for the averaged iterate in both model-based and model-free {LQ RL}, and provide an explicit expression for the corresponding asymptotic variance, along with the exact asymptotic behavior of the objective function loss. Notably, the analysis of objective loss is motivated by the observation that, in contrast to many {ML} settings, where the objective loss often lacks direct practical interpretation, the loss in {LQ RL} can carry concrete significance in specific application domains. For instance, in portfolio liquidation, the loss corresponds to expected returns, while in clinical applications such as heart disease management, it reflects a patient's risk score. Moreover, from a practical standpoint, one may be primarily concerned with the anticipated level of risk or cost, and only secondarily with the precise policy used to achieve it.

\begin{theorem}[Asymptotic normality for model-based {LQ RL}]\label{prop:clt1}
	Under the  Assumptions \ref{assu1}-\ref{assu5}, we have
	\begin{equation*}
		\frac1{\sqrt n}\sum_{\ell=1}^n\left(\operatorname{vec}\left\{K_t^\ell\right\}-\operatorname{vec}\left\{K_t^{\star}\right\}\right) \xrightarrow{d} \mathcal{N}\left(0, G^{-1}_t S_t G^{-1}_t\right),
	\end{equation*}
	where $G_t=\mathbb{E}\left[x_tx_t^\top\right]\otimes\left(R_t+B^\top P^{\star}_{t+1}B\right)$, $S_t=\mathbb{E}\left[\operatorname{vec}\left\{E_t\left(K_t^{\star}\right)x_t x_t^\top\right\}\operatorname{vec}\left\{E_t\left(K_t^{\star}\right)x_t x_t^\top\right\}^\top\right]$, and $x_t \sim \mathcal{P}_t$.
\end{theorem}


\begin{corollary}[Asymptotics of objective function loss] \label{prop:loss_clt1}
	Under the assumptions as listed in Theorem \ref{prop:clt1}, the objective loss follows a mixed Chi-squared distribution asymptotically,
	\begin{align*}
		2n\left(C\left(\boldsymbol{\bar{K}}\right) - C\left(\boldsymbol{K}^\star\right)\right) \xrightarrow{d}  Z^\top H Z \sim  \sum_{t = 1}^{T} \lambda_t \chi^2_t(1),
	\end{align*}
	where $Z = \sqrt{n} \left(\operatorname{vec}\left\{\boldsymbol{\bar{K}}\right\} - \operatorname{vec}\left\{\boldsymbol{K}^\star\right\}\right) \xrightarrow{d} \mathcal{N}(0, \boldsymbol{G}^{-1}\boldsymbol{S}\boldsymbol{G}^{-1})$, $H = \nabla^2 C(\boldsymbol{K}^\star)$ is the block-diagonal Hessian with $\nabla^2_{t} C(\boldsymbol{K}^\star) = 2 G_t \otimes \left(R_t + B^\top P^{\star}_{t+1} B\right)$, and $\lambda_t$ is the eigenvalue of $G^{1/2}HG^{1/2}$. We note that using the same argument as in the proof of Theorem \ref{prop:loss_clt1}, one can immediately derive CLT results of objective loss for the model-free setting.
\end{corollary}


\begin{theorem}[Asymptotic normality for model-free {LQ RL}]\label{prop:clt3}
	Under the Assumptions \ref{assu1} to \ref{assu5}, we have
	\begin{eqnarray*}
		\frac1{\sqrt n}\sum_{\ell=1}^n\left(\vec\left\{\tilde K_t^\ell\right\}-\vec\left\{\tilde K_t^{\star}\right\}\right)\xrightarrow{d} \mathcal{N}(0, \tilde G^{-1}_t \tilde S_t \tilde G^{-1}_t),
	\end{eqnarray*}
	where 
	\begin{eqnarray*}
		\tilde G_t &=& \E\left[x_tx_t^\top\right]\otimes\left(R_t+B^\top P^{\star}_{t+1}B\right), \\
		\tilde S_t &=& \frac1{m^2}\sum_{i,j=1}^m\E\left[\vec\left\{\left(c\left(\boldsymbol{K}^{\star}+ U_t^i\right)-c\left(\boldsymbol{K}^{\star}\right)\right) U_t^i\right\} \vec\left\{\left(c\left(\boldsymbol {K}^{\star}+ U_t^j\right)-c\left(\boldsymbol {K}^{\star}\right)\right) U_t^j\right\}^\top\right],
	\end{eqnarray*}
	and $x_t \sim \mathcal{P}_t$.    
\end{theorem}


\begin{remark}
	The proof of {CLT} results is inspired by the results in \cite{polyak1992acceleration}, utilizing the vectorization technique. To establish the asymptotic consistency of the proposed bootstrap algorithm, we compare the tail probability of the true distribution and the bootstrap distribution, where we construct respective expansions that satisfy asymptotic normality. Using the constructed Gaussian process as a bridge for the model-based setting, we then analyze the Kolmogorov distance between the true distribution and the bootstrap distribution. On the other hand, for the model-free setting, we uniformly perturb the cost function and estimate its gradient using empirical observational data. We then employ similar techniques in the model-based case to conduct the analysis.
\end{remark}

Theorems \ref{prop:clt1} and \ref{prop:clt3} provide the statistical foundation for the inference procedure of ${K}^{\star}_t$. Given the explicit form of covariance as in Theorem \ref{prop:clt1}, the asymptotic covariance may be estimated using a plug-in estimator by using $\bar{K}^n_t$ and the sample covariance of the state. However, in the context of noisy {LQ RL} systems, the randomness of the state comes from the noise term and the initial state. The joint distribution of these two terms is rarely known in the real world. Our bootstrapping algorithm provides an efficient way to construct a confidence interval for the optimal policy $K^{\star}_t$. The asymptotic validation of the bootstrapping algorithm is established in the following section.

\subsection{Asymptotic validation of bootstrap algorithm}

In this section, we focus on developing conditional convergence for bootstrap estimation, which provides the theoretical guarantee of the asymptotic validity of the constructed confidence interval and other inferential tasks. 

\begin{theorem}[Asymptotic expansion of exact {PG} iterates]\label{lem:clt2}
	Under the Assumption \ref{assu1}-\ref{assu5}, we have
	\begin{equation*}\label{e:dec}
		\frac1{\sqrt n}\sum_{\ell=1}^n\left(\operatorname{vec}\left\{K_t^\ell\right\}-\operatorname{vec}\left\{K_t^{\star}\right\}\right)= \frac1{\sqrt n} G_t^{-1}\sum_{\ell=1}^n{\xi}_t^\ell(0) + o_p(1),
	\end{equation*}	
	where ${\xi}_t^\ell(0) = -2W_\ell \operatorname{vec}\left\{\left(\left(R_t+B^{\top} P_{t+1} B\right)K_t^{\star}-B^{\top} P^{\boldsymbol{K}}_{t+1} A \right) x_t^\ell \left(x_t^\ell\right)^\top\right\}$.
\end{theorem}			


\begin{theorem}[Asymptotic expansion of zeroth-order {PG}]\label{lem:clt4}
	Under the Assumptions \ref{assu1}-\ref{assu5}, we have
	\begin{eqnarray}\label{e:dec2}
		\frac1{\sqrt n}\sum_{\ell=1}^n\left(\vec\left\{\tilde K_t^\ell\right\}-\vec\left\{ K_t^{\star}\right\}\right)= \frac1{\sqrt n} \tilde G_t^{-1}\sum_{\ell=1}^n{\tilde \xi}_t^\ell(0) + o_p(1),
	\end{eqnarray}	
	where ${\tilde \xi}_t^\ell(0)= -\frac{W_\ell}{m}\sum_{i=1}^{m}\vec\left\{\left(c\left(\boldsymbol{K}^{\star}+ U_t^i\right)-c\left(\boldsymbol {K}^{\star}\right)\right)U_t^i\right\}.$
\end{theorem}


Based on the above asymptotic expansion, we can see that $\sqrt{n} \left(\vec\left\{\bar{K}_t^n\right\} - \vec\left\{K_t^{\star}\right\}\right)$ has an asymptotic Gaussian distribution. In addition, the distribution can be expected to enjoy asymptotic consistency, i.e., the perturbed empirical distribution of the proposed bootstrap estimator should converge to the true underlying distribution at the rate of $\sqrt{n}$. We formalize this conclusion by the following theorem.

\begin{theorem}[Asymptotic Validation of Bootstrap Algorithm (Model-based)]\label{thm:model base}
	Under the Assumptions \ref{assu1}-\ref{assu5}, and $n \rightarrow \infty$, we have 
	{\small\begin{equation*}
			\sup_{v \in \mathbb{R}^{dp}}\left|\mathbb{P}^\star\left(\frac1{\sqrt n}\sum_{\ell=1}^n\left(\operatorname{vec}\left\{K_t^{\ell, \prime}\right\}-\operatorname{vec}\left\{K_t^{\ell}\right\}\right)\le v\right)- \\
			\mathbb{P}\left(\frac1{\sqrt n}\sum_{\ell=1}^n\left(\operatorname{vec}\left\{K_t^\ell\right\}-\operatorname{vec}\left\{K_t^{\star}\right\}\right)\le v \right)\right| 
			\to 0
	\end{equation*}}
	in probability.
\end{theorem}


\begin{theorem}[Asymptotic Validation of Bootstrap Algorithm (Model-free)]\label{thm:model free}
	Under the Assumptions \ref{assu1}-\ref{assu4}, and $n \rightarrow \infty$, we have 
	\begin{align*}
		\sup_{v \in \mathbb{R}^{dp}}\left|\mathbb{P}^\star\left(\frac1{\sqrt n}\sum_{\ell=1}^n\left(\vec\left\{\tilde K_t^{\ell, \prime}\right\}-\vec\left\{\tilde K_t^\ell\right\}\right)\le v\right)- \mathbb{P}\left(\frac1{\sqrt n}\sum_{\ell=1}^n\left(\vec\{\tilde K_t^\ell\}-\vec\{K_t^{\star}\}\right)\le v \right)\right| \\ 
		\to 0	
	\end{align*}
	in probability.
\end{theorem}


Theorem \ref{lem:clt2} and \ref{thm:model base} verify the validity and suitability of our proposed bootstrapping algorithm to conduct statistical inference. Theorem \ref{thm:model base} establishes the distributional convergence of our bootstrap algorithm in terms of the Kolmogorov metric, enabling us to construct a confidence interval for the optimal policy $K^{\star}_t$ that utilizes the bootstrap samples $\left\{\bar{K}_{t, (b)}^{n, \prime}\right\}_{b=1}^\mathcal{B}$. In practice, an inherent Monte Carlo error arises due to resampling from a finite number of bootstrapped samples $\mathcal{B}$. However, numerical experiments show that this error becomes negligible as $\mathcal{B}$ grows sufficiently large.

We emphasize that the results established in Theorems \ref{lem:clt2} and \ref{thm:model base} are primarily important among the asymptotic ones of the bootstrap estimates, which play a crucial role in validating (asymptotic) confidence intervals and conducting hypothesis testing. However, the natural question regarding the exact rate of convergence in the relevant consistency remains insufficiently addressed in the literature. To this end, we establish a non-asymptotic analysis in the next section to justify confidence intervals for the underlying true policy in LQ RL using the multiplier bootstrap.

\subsection{Non-asymptotic validation of Multiplier Bootstrap}

In this section, we establish an exact finite-sample derivation for the accuracy of the bootstrap approximation of the distribution of Polyak-Ruppert averaging iterates with a polynomially decreasing step size, leveraging Bernstein-type concentration inequalities. The obtained results provide significant insight into the fact that the quantiles of the exact distribution of $\sqrt{n}\left(\vec\{\bar K_t^n\} - \vec\{ K_t^{\star}\}\right)$ can be approximated at a best rate of $n^{-1/4}$, where $n$ is the number of agents used in the procedure,  provided that $n$ is sufficiently large. To conduct a non-asymptotic analysis of the bootstrap procedure, we impose the following additional essential assumption.

\begin{assumption}\label{assu6}
	We denote     
	\begin{eqnarray*}
		\epsilon_t^n&=&2\left(x_t^n\left(x_t^n\right)^\top-\E\left[x_t^n\left(x_t^n\right)^\top\right]\right)\otimes \left(R_t+B^\top P_{t+1}B\right)\vec\{K_t^\star\}\\
		&&-2\vec\left\{B^\top P_{t+1}A\left(x_t^n\left(x_t^n\right)^\top-\E\left[x_t^n\left(x_t^n\right)^\top\right]\right)\right\}, 
	\end{eqnarray*} 
	where we assume that its covariance matrix exists and has the smallest eigenvalue $\lambda_\epsilon>0$.    
\end{assumption}

\begin{assumption}\label{assu7}
	The step size is chosen as $\eta_\ell=c_0\ell^{-\frac12}$ for a small enough constant $c_0$ and $\frac{\sqrt{n}}{(1+\log(n))\log (n)}\ge C_{0}$, where the constant $C_0$ is a strict positive number depends on the LQ RL parameters. Moreover, we assume that
	\begin{eqnarray*}
		\lambda_{\epsilon}\ge 8C_\lambda\sqrt{\frac{\left|\operatorname{cov}\left(\epsilon_t^n\right) \right|\log (n)}{n}} +\frac{8\left(\left|\operatorname{cov}\left(\epsilon_t^n\right)\right|+C_\lambda^2\right)\log (n)}{n}.
	\end{eqnarray*}
\end{assumption}

\begin{theorem} [Non-asymptotic Validation of Bootstrap Algorithm (Model-based)] \label{them:nonasy}
	Under extra Assumptions \ref{assu6} \& \ref{assu7}, we have with probability at least $1-6/n$,
	\begin{eqnarray*} 
		&&\sup_{v \in \mathbb{R}^{dp}}\left|\mathbb{P}^{\star}\left(\sqrt{n}\left(\operatorname{vec}\left\{\bar{K}_t^{n, \prime}\right\}-\operatorname{vec}\left\{\bar{K}_t^n\right\}\right) \le v\right)- \mathbb{P}\left(\sqrt{n}\left(\operatorname{vec}\left\{\bar{K}_t^n\right\}-\operatorname{vec}\left\{K_t^{\star}\right\}\right) \le v \right)\right|\\
		&\lesssim & \frac{\kappa^2 \left(C_A^4 \vee 1 \right)\left(1 + C_\epsilon^2\right) \log(n)}{\lambda_\epsilon n^{1/4}} 
		+ \frac{\sqrt{dp}}{\sqrt{n}} \left( \frac{C_\epsilon^3}{\lambda_\epsilon^{3/2}} + \kappa C_\epsilon \frac{\sqrt{\log(n)}}{\sqrt{\lambda_{\epsilon}}} + \frac{\kappa \left(1 + C_\epsilon^2 / \lambda_{\epsilon}\right) \log (n)}{\sqrt{n}} \right)\\ 
		&& + \frac{C_{n,\lambda_\epsilon,C_A} e^{-C\sqrt{n}}}{\lambda_{\epsilon}} \left|\operatorname{vec}\left\{K_t^0\right\} - \operatorname{vec}\left\{K_t^\star\right\}\right|,
	\end{eqnarray*}
	where $\kappa, C_A$ are positive constants that depend on the LQ RL, $C_{n,\lambda_\epsilon, C_A}$ is a constant that depends on $n$ polynomially.
\end{theorem}


We emphasize that the non-asymptotic results of Theorem \ref{them:nonasy} constitute a key step in establishing the validity of the bootstrap procedure.  Theorem \ref{them:nonasy} justifies that the best achievable rate of convergence for the bootstrap procedure has an error order of $n^{-1/4}$, which unfortunately does not attain the statistically optimal rate of $n^{-1/2}$. This result is valuable in filling the gap in understanding the bootstrap procedure, specifically regarding its sample efficiency. Under finite-sample settings, it provides guidance for practical implementation by highlighting that, when the sample size is limited, the bootstrap may lack coverage guarantees, and consequently, the resulting confidence intervals should be interpreted with caution.

\begin{corollary} 
	
	We consider the problem of estimating the $\alpha$-quantile, for some $\alpha \in(0,1)$, associated with a matrix $\boldsymbol{D} \in$ $\mathbb{R}^{d p \times d p}$, formally defined as the value for time point $t$. 
	\begin{equation*}
		\mu_\alpha=\inf \left\{\mu>0: \mathbb{P}\left(\sqrt{n}\left\|\boldsymbol{D}\left(\operatorname{vec}\left\{\bar{K}_t^n\right\} - \operatorname{vec}\left\{K_t^{\star}\right\}\right)\right\| \geq \mu\right) \leq \alpha\right\}.
	\end{equation*}
	Define its counterpart in the Bootstrap world, that is, the quantity$t_\alpha^{\prime}$, as 
	\begin{equation*}
		\mu_\alpha^{\prime}=\inf \left\{\mu>0: \mathbb{P}\left(\sqrt{n}\left\|\boldsymbol{D}\left(\operatorname{vec}\left\{\bar{K}_t^{n, \prime}\right\} - \operatorname{vec}\left\{\bar{K}_t^{n}\right\}\right)\right\| \geq \mu\right) \leq \alpha\right\}. 
	\end{equation*}
	We clarify that $\mu_\alpha^{\prime}$ is defined w.r.t the bootstrap measure; this bootstrap critical value $\mu_\alpha^{\prime}$ is then used within the bootstrap empirical distribution to construct the confidence set.
	\begin{equation*}
		\mathcal{E}(\alpha)=\left\{\operatorname{vec}\{K_t\} \in \mathbb{R}^{dp}: \sqrt{n}\left\|\boldsymbol{D}\left(\operatorname{vec}\{K_t\} - \operatorname{vec}\left\{\bar{K}_t^n\right\}\right)\right\| \leq \mu_\alpha^{\prime}\right\}
	\end{equation*}
	It states that the non-asymptotic error bound provides a rigorous assessment of the probability coverage guarantee that the confidence set includes the true parameter $\operatorname{vec}\left\{K^{\star}_t\right\}$, that is, 
	\begin{equation*}
		\mathbb{P}\left(\operatorname{vec}\left\{K^{\star}\right\} \notin \mathcal{E}(\alpha)\right)=\mathbb{P}\left(\sqrt{n}\left\|\boldsymbol{D}\left(\operatorname{vec}\left\{\bar{K}_t^n\right\} - \operatorname{vec}\left\{K_t^{\star}\right\}\right)\right\|>\mu_\alpha^{\prime}\right) \approx \alpha,
	\end{equation*}
	which deviates from the nominal level by an error term on the order of $n^{-1 / 4}$.
\end{corollary}

\begin{remark} We note that the finite-sample analysis conducted here builds upon the non-asymptotic rates for martingale properties of SGD established in \cite{anastasiou2019normal}. Aside from the notable exception of \cite{samsonov2024gaussian}, which primarily focuses on linear stochastic approximation (LSE), we are not aware of any existing work that investigates the non-asymptotic properties of the bootstrap for Polyak-Ruppert averaging iterates. Nonetheless, the complexity of LQ RL introduces additional challenges beyond those encountered in LSE when addressing this question. Importantly, such an argument has been missing in earlier works constructing confidence intervals using the bootstrap for statistical inference in reinforcement learning (e.g., \cite{shi2023dynamic, ramprasad2023online}).
	
	We expect that the proof technique developed for Theorem \ref{them:nonasy} can be extended to provide similar non-asymptotic validity results for various model-free LQ RL approaches. However, we emphasize that generalizing Theorem \ref{them:nonasy} to model-free settings—where gradients are approximated via zeroth-order gradient estimators—is substantially more complex, as it critically depends on finite-sample approximation error rates of these estimators. Consequently, a fundamental challenge arises from the lack of established non-asymptotic error bounds for zeroth-order gradient estimation. The proof techniques employed here are not directly applicable, and developing an appropriate non-asymptotic analysis for zeroth-order gradient estimators constitutes a distinct and challenging research direction. Although this topic is important, our primary focus is on conducting statistical inference rather than investigating the theoretical properties of learning algorithms, which we leave as a direction for future research.
\end{remark}

\section{Numerical experiments} \label{sec:ne}

In this section, we evaluate the performance of the proposed method for both model-based and model-free {LQ RL} and consider the settings as follows, 
\begin{equation*}
	A=\left(\begin{array}{cccc}
		0.5 & 0.05 & 0.1 & 0.2 \\ 
		0 & 0.2 & 0.3 & 0.1\\
		0.06 & 0.1 & 0.2 & 0.4\\
		0.05 & 0.2 & 0.15 & 0.1
	\end{array}\right), 
	B=\left(\begin{array}{cc}
		-0.05 & -0.01\\ 
		-0.005 & -0.01 \\
		-1 & -0.01\\
		-0.01 & 0.9\end{array}\right), 
	R_t=\left(\begin{array}{cc}
		0.4 & -0.25\\ 
		-0.25 & 0.7 
	\end{array}\right), 
\end{equation*} 

\begin{equation*}
	Q_T=\left(\begin{array}{cccc}
		1 & 0.2 & -0.005 & 0.015\\
		0.2 & 1.1 & 0.15 & 0\\
		-0.05 & 0.15 & 0.9 & -0.08\\
		0.015 & 0 & -0.08 & 0.88
	\end{array}\right),
	W=\left(\begin{array}{cccc}
		0.1 & 0 & 0 & 0\\
		0 & 0.5 & 0 & 0\\
		0 & 0 & 0.2 & 0\\
		0 & 0 & 0 & 0.3
	\end{array}\right).
\end{equation*} 
We consider the case that $T = 10$ and $Q_t = Q_T, \forall t$, $x_0=(x^1_0,x^2_0,x^3_0,x^4_0)^\top$ and $x^i_0$ are independent, and $x^1_0,x^2_0,x^3_0,$ and $x^4_0$ are sampled from $\mathcal{N}(5,0.1),$ $\mathcal{N}(2,0.3),$ $\mathcal{N}(8,1),$ $\mathcal{N}(5,0.5)$, respectively. We use the line search method to choose the step size. We mainly focus on the model-based setting to explore asymptotic consistency and validation (see Appendix 
for the results of the model-free cases).

\subsection{Convergence of averaging iterates} 

We first demonstrate the performance of the Polyak-Ruppert averaging iterates of policy gradient for the setting with $m = 20$ with initial policy $\{K^0\}_{ij}=0.4$ for all $i, j$. 

To show the convergence of policy $\boldsymbol{\bar{K}}^n$, we consider \textit{bias} $\|\operatorname{vec}(\boldsymbol{\bar{K}}^n) - \operatorname{vec}(\boldsymbol{K}^{\star})\|_1$, \textit{cost of execution} $C(\boldsymbol{K}^n)$, and \textit{normalized error} $\frac{C(\boldsymbol{\bar{K}}^n)-C\left(\boldsymbol{K}^{\star}_t\right)}{C\left(\boldsymbol{K}^{\star}\right)}$, where $K^{\star}$ is the optimal policy. 

\begin{figure}[htp!]
	\centering
	\includegraphics[scale = 0.17]{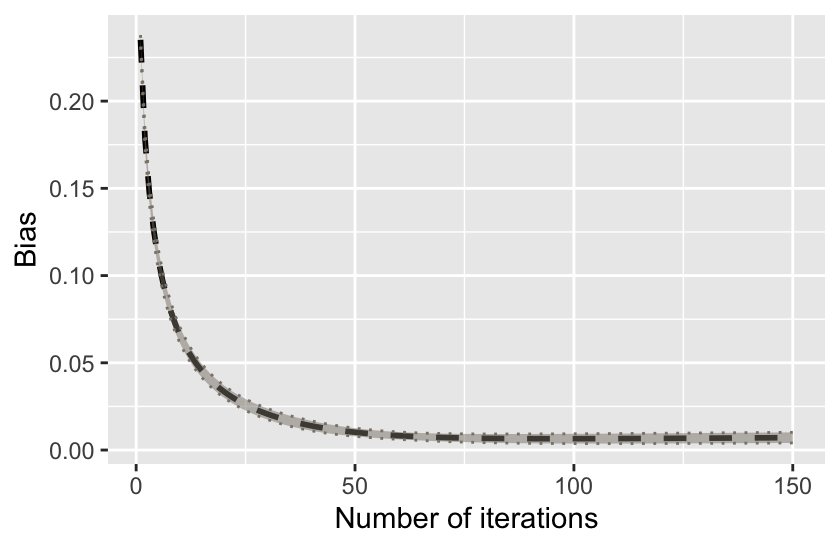}
	\includegraphics[scale = 0.17]{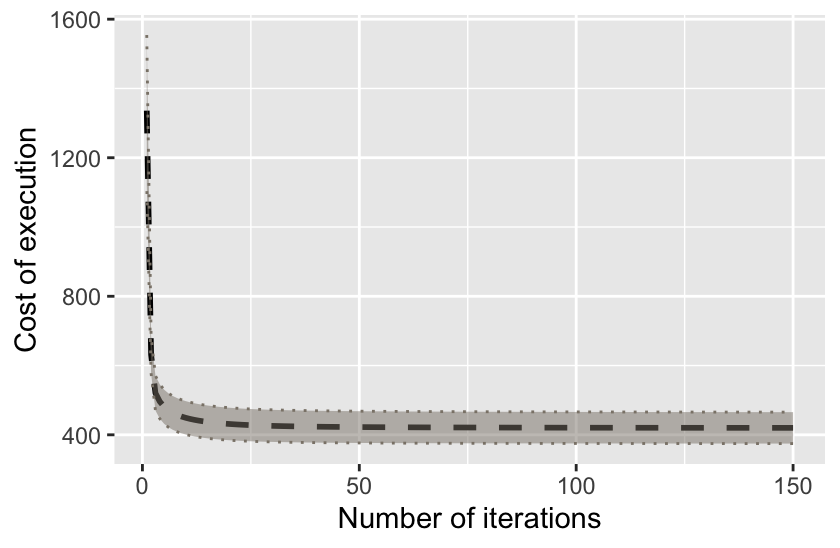}
	\includegraphics[scale = 0.17]{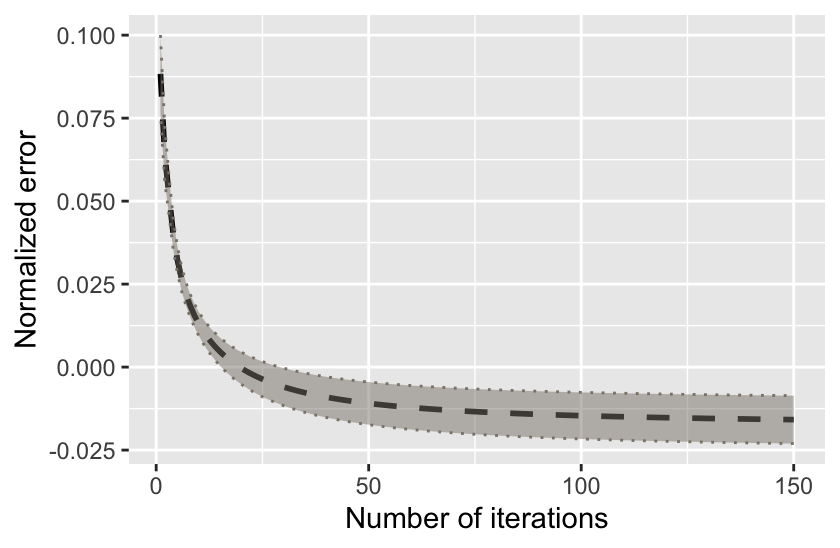}
	\caption{\textit{Model-based: Bias, loss function, and normalized error for $\boldsymbol{K}^n_t$ when $M = 20$ with $100$ replications. The shadow area on the right stands for the variability of the normalized error.}}
	\label{covergence}
\end{figure}

Figure \ref{covergence} summarizes the results of multiple agents under the model-based setting ($m = 20$). It illustrates that the exact {PG} averaging iterates converge with a reasonable level of accuracy guarantee within $100$ iterations (the bias is less than $10^{-2}$). Performance is stable with relatively minor fluctuations (see the shadow area in Figure \ref{covergence}) in the $100$ replications. Notably, {PG}'s convergence is sensitive to the initial value of $\boldsymbol{K}^0$, the initial state, and the stepsize. The learning trajectories have different shapes, with various initial states and step sizes.

\begin{figure}[t!]
	\centering
	\includegraphics[scale = 0.3]{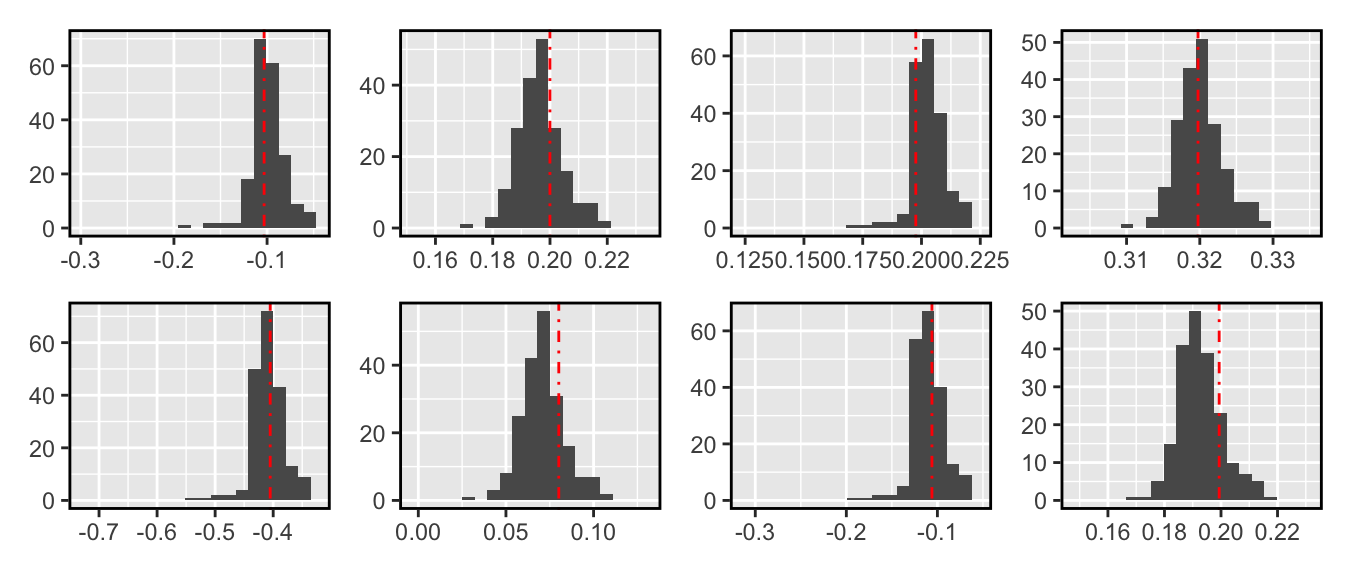}
	\caption{\textit{Model-based: The histogram for the bootstrap estimates for all dimensions of $K$ at $t=1$, {i.e.}, $\{{K}_1\}_{ij}$, the dot-dash line is the true value.}}
	\label{hist}
\end{figure}

\begin{figure}[t!]
	\centering
	\includegraphics[scale = 0.13]{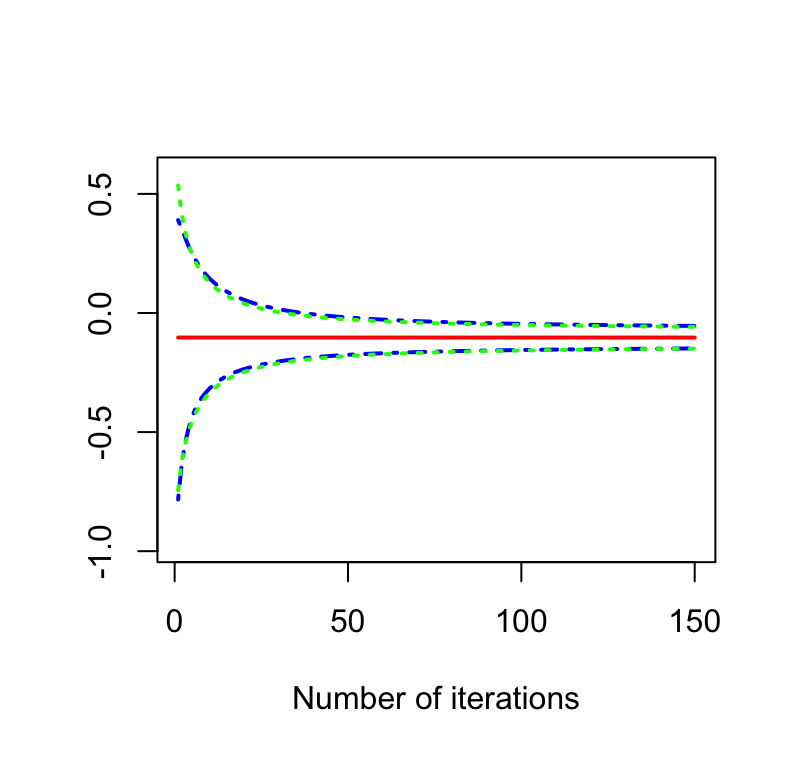}
	\includegraphics[scale = 0.13]{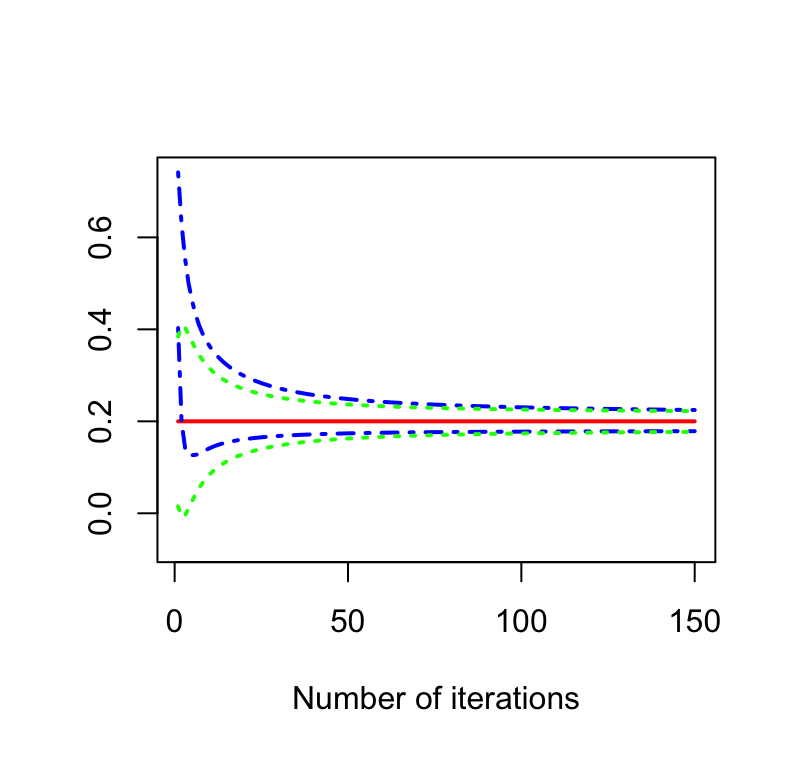}
	\includegraphics[scale = 0.13]{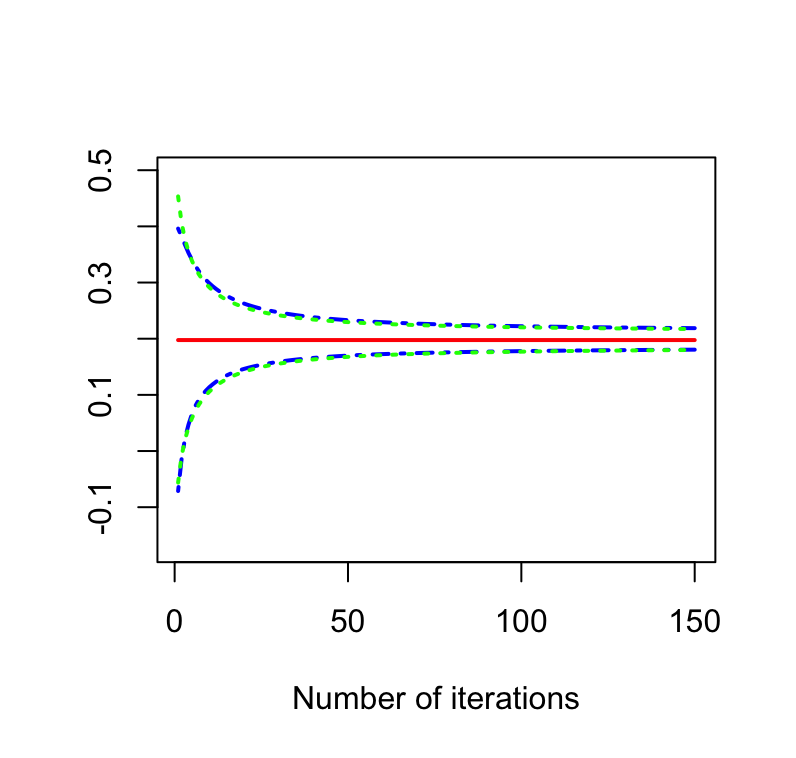}
	\includegraphics[scale = 0.13]{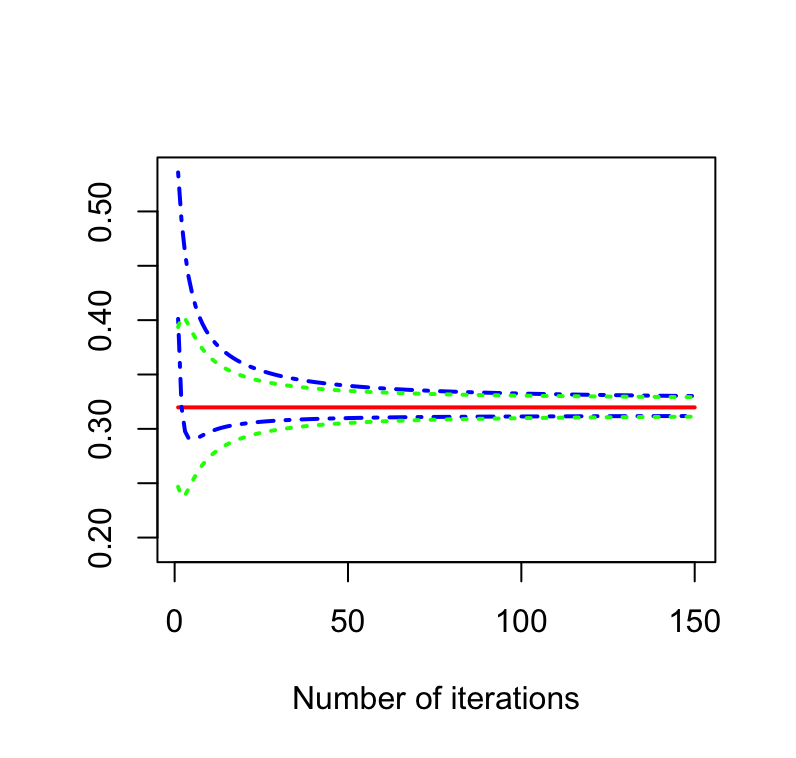}\\
	\includegraphics[scale = 0.13]{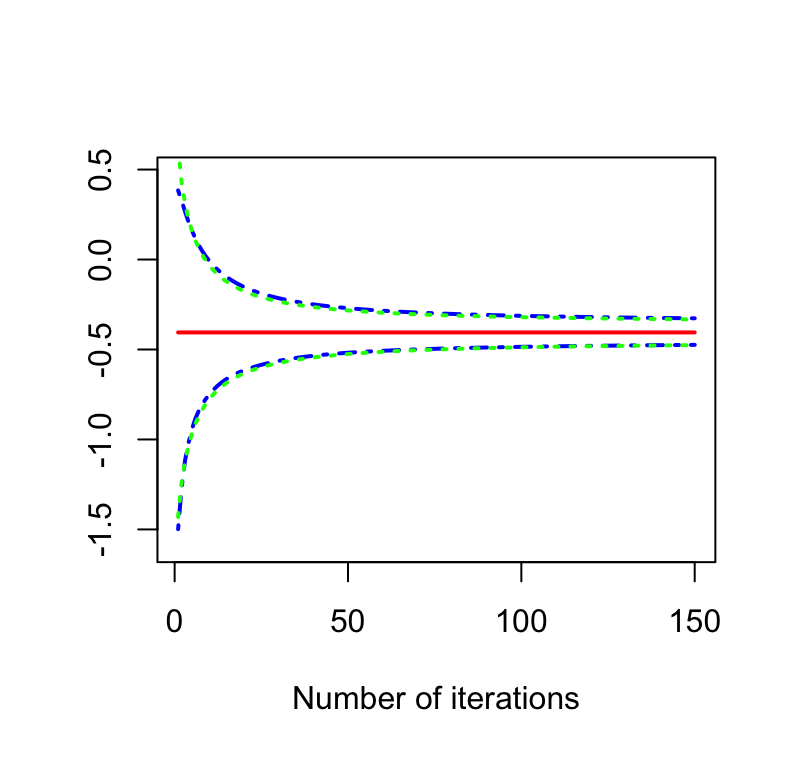}
	\includegraphics[scale = 0.13]{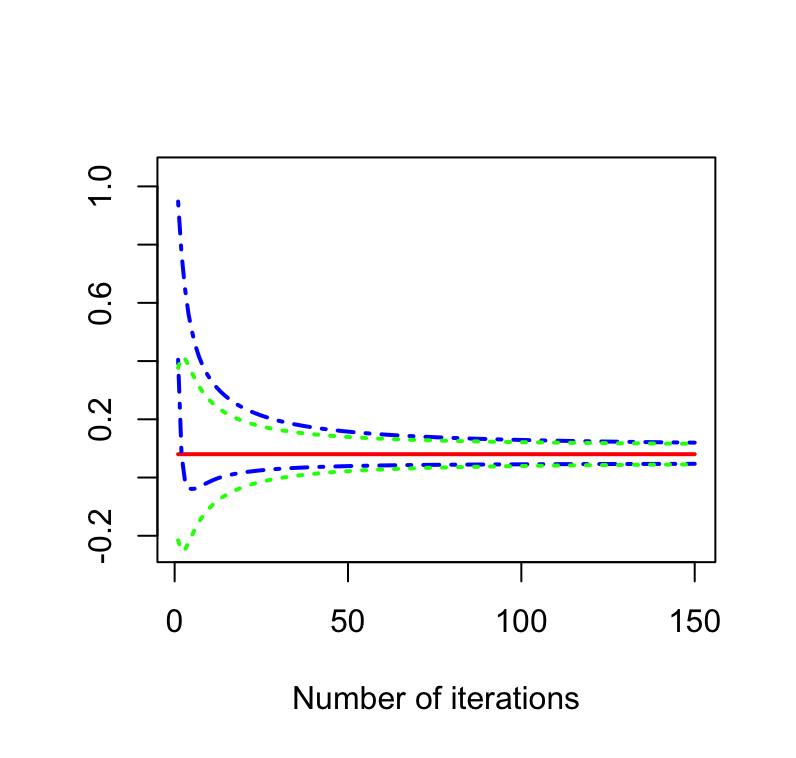}
	\includegraphics[scale = 0.13]{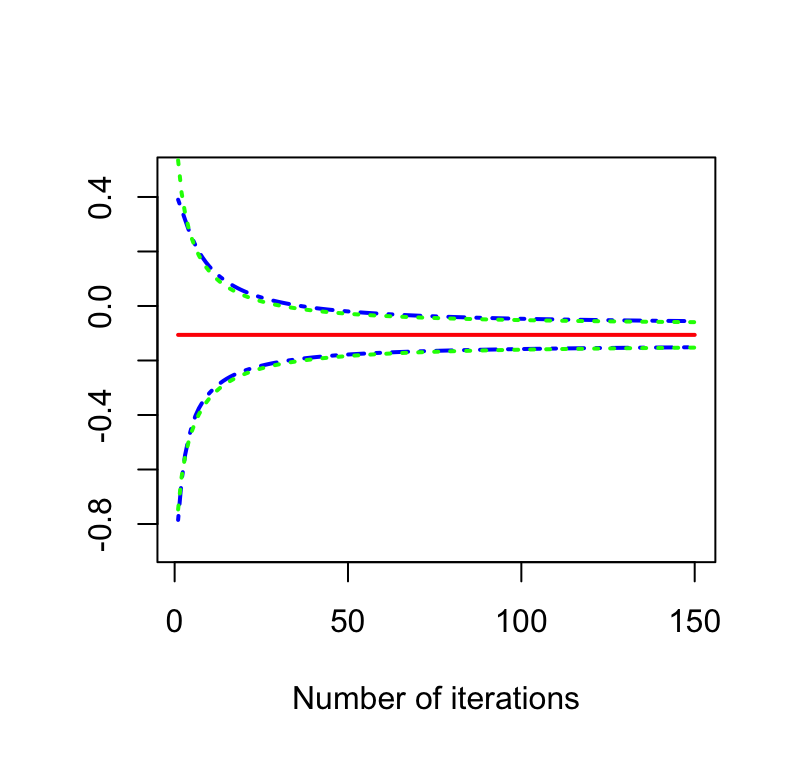}
	\includegraphics[scale = 0.13]{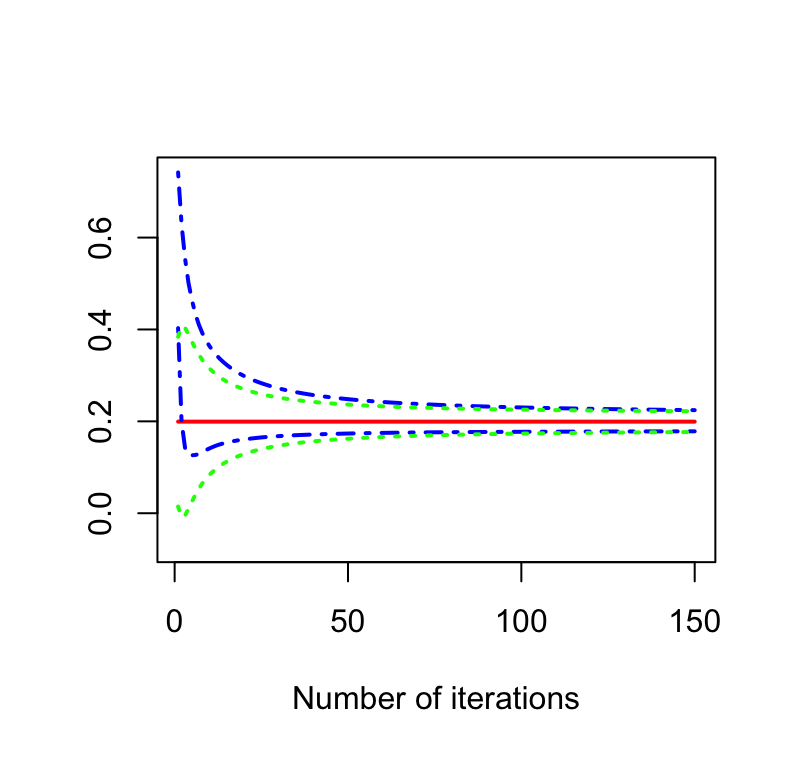}\\
	\caption{\textit{Model-based: The confidence band and true value of all the dimensions of $K$ at $t=1$ with $100$ replications, where $M = 20$, $B=200$. From left to right, the top three figures: $\{K_{1}\}_{11}, \{K_{1}\}_{21},\{K_{1}\}_{12}$; Middle three figures: $\{K_{1}\}_{22}, \{K_{1}\}_{31},\{K_{1}\}_{32}$; Bottom two: $\{K_{1}\}_{41}, \{K_{1}\}_{42}$.}}
	\label{ci_modelbased}
\end{figure}

\subsection{Online bootstrapping inference for {PG} estimator} \label{sec:validaton}

In this section, we evaluate the performance of the proposed bootstrapping procedure to construct quantile-based and standard-error-based confidence intervals for the optimal policy $\boldsymbol{K}^{\star}$. To perform the bootstrapping procedure, we consider $B = 200$ and $W_{n, (b)}$ to be generated from an exponential distribution with mean and variance of $1$. The results are summarized in Figures \ref{hist} - \ref{ci_modelbased}, where the solid line indicates the true value and, for Figure \ref{ci_modelbased}, the dot-dash line shows the confidence band obtained based on the quantile value of $0.025$th and $0.975$th bootstrap estimates, and the dotted line indicates the confidence band obtained based on the standard variance of estimates.

We can see from Figure \ref{hist} that the histogram is close to normal sharpness, and the true value $K^\star$ is almost in the middle of the histogram in all dimensions. It provides empirical evidence for the asymptotic normality of the policy gradient estimator established in Section \ref{sec:mainresult}. Figure \ref{ci_modelbased} shows that the confidence intervals derived from the bootstrapping policy gradient descent method successfully cover the true value. It further confirms the distributional consistency of the bootstrapping algorithm, {i.e.}, as the number of iterations increases, the length of the confidence interval constructed decreases, and the empirical distribution of the perturbed {PG} estimates converges to a normal distribution (the quantile-based confidence interval constructed coincides with the standard error-based confidence interval constructed based on normal application).  

\subsection{Sensitivity analysis}

We further investigate the impact of the weighting schemes and the underlying dynamics on the empirical performance of the proposed bootstrapping procedure. We denote $RW-Q$ for the quantile-based confidence interval, where the proposed random weighting procedure is used to obtain quantiles of $2.5\%$ and $97.5\%$, and $RW-\sigma$ for the standard-error-based confidence interval, where the proposed random weighting procedure is used to estimate its standard error. There are two types of random weighting schemes under consideration: (a). Weights generated from a normal distribution, $W_{n, (b)} \sim \mathcal{N}(1, 1)$ (an absolute value to be used if generated $W_{n, (b)}$ is negative); (b). Weights generated from an exponential distribution, $W_{n, (b)} \sim {Exp}(1)$, under two types of stochastic linear dynamics, (a). Random noise, $\omega_t$, generated from a normal distribution; (b). Random noise, $\omega_t$, generated from a $t$-distribution.

\begin{table}[ht!]
	\centering 
	\resizebox{\columnwidth}{!}{\begin{tabular}{c c c cccccccc} 
			\hline\hline Number of agent, $m$ & Approach & Weight scheme & $K_{1,1}$ & $K_{1,2}$ & $K_{1,3}$ & $K_{1,4}$ & $K_{2,1}$ & $K_{2,2}$ & $K_{2,3}$ & $K_{2,4}$ \\
			\hline
			\multirow[t]{4}{*}{$$
				m=10
				$$
			} & $R W-Q$ & ${Exp}(1)$ & 0.869 & 0.890 & 0.907 & 0.918 & \textbf{0.950} & 0.942 & 0.903 & 0.909 \\
			& & $\mathcal{N}(1,1)$ & 0.794 & 0.840 & 0.753 & 0.829 & 0.873 & 0.882 & 0.781 & 0.815 \\
			& $R W-\sigma$ & ${Exp}(1)$ & 0.868 & 0.892 & 0.899 & 0.919 & 0.942 & 0.942 & 0.905 & 0.917 \\
			& & $\mathcal{N}(1,1)$ & 0.865 & 0.893 & 0.815 & 0.862 & 0.917 & 0.921 & 0.827 & 0.858 \\
			\hline \multirow[t]{4}{*}{$$
				m=20
				$$
			} & $R W-Q$ & ${Exp}(1)$ & \textbf{0.956} & \textbf{0.966} & 0.937 & 0.939 & \textbf{0.984} & \textbf{0.967} & \textbf{0.945} & \textbf{0.946} \\
			& & $\mathcal{N}(1,1)$ & 0.883 & 0.916 & 0.848 & 0.891 & 0.938 & 0.935 & 0.880 & 0.904 \\
			& $R W-\sigma$ & ${Exp}(1)$ & \textbf{0.951} & \textbf{0.960} & 0.940 & \textbf{0.951} & \textbf{0.979} & \textbf{0.970} & \textbf{0.953} & \textbf{0.952} \\
			& & $\mathcal{N}(1,1)$ & 0.927 & \textbf{0.946} & 0.865 & 0.919 & \textbf{0.962} & \textbf{0.952} & 0.917 & 0.930 \\
			\hline\hline
	\end{tabular}}
	\caption{\textit{Model-based: Average coverage probabilities with the coverage level being $0.95$ with normal distributed system error $\omega_t$. }}
	\label{tab:t_dis}
\end{table}

\begin{table}[ht!]
	\centering 
	\resizebox{\columnwidth}{!}{\begin{tabular}{c c c cccccccc} 
			\hline\hline Number of agent, $m$ & Approach & Weight scheme & $K_{1,1}$ & $K_{1,2}$ & $K_{1,3}$ & $K_{1,4}$ & $K_{2,1}$ & $K_{2,2}$ & $K_{2,3}$ & $K_{2,4}$ \\
			\hline \multirow[t]{4}{*}{$$
				m=20
				$$
			} & $R W-Q$ & ${Exp}(1)$ & 0.876 & 0.887 & 0.854 & 0.892 & 0.914 & 0.916 & 0.900 & 0.901 \\
			& & $\mathcal{N}(1,1)$ & 0.867 & 0.859 & 0.728 & 0.766 & 0.852 & 0.866 & 0.829 & 0.834 \\
			& $R W-\sigma$ & ${Exp}(1)$ & 0.891 & 0.901 & 0.857 & 0.897 & 0.919 & 0.922 & 0.902 & 0.901 \\
			& & $\mathcal{N}(1,1)$ & 0.858 & 0.874 & 0.752 & 0.813 & 0.874 & 0.882 & 0.867 & 0.874 \\
			\hline \multirow[t]{4}{*}{$$
				m=50
				$$
			}&  $R W-Q$ & ${Exp}(1)$ & 0.928 & 0.910 & 0.924 & 0.907 & \textbf{0.965} & 0.933 & \textbf{0.946} & 0.926 \\
			& & $\mathcal{N}(1,1)$ & 0.874 & 0.867 & 0.837 & 0.814 & 0.900 & 0.900 & 0.892 & 0.891 \\     
			& $R W-\sigma$ & ${Exp}(1)$ & 0.938 & 0.914 & 0.919 & 0.914 & \textbf{0.958} & 0.926 & 0.937 & 0.924 \\
			& & $\mathcal{N}(1,1)$ & 0.909 & 0.901 & 0.857 & 0.872 & 0.942 & 0.926 & 0.914 & 0.914 \\
			\hline \multirow[t]{4}{*}{$$
				m=100
				$$
			} & $R W-Q$ & ${Exp}(1)$ & \textbf{0.981} & \textbf{0.948} & \textbf{0.970} & 0.936 & \textbf{0.995} & \textbf{0.960} & \textbf{0.956} & \textbf{0.956} \\
			& & $\mathcal{N}(1,1)$ & 0.893 & 0.888 & 0.871 & 0.826 & 0.927 & 0.924 & 0.896 & 0.917 \\
			& $R W-\sigma$ & ${Exp}(1)$ & \textbf{0.988} & \textbf{0.951} & \textbf{0.966} & 0.939 & \textbf{0.994} & \textbf{0.957} & \textbf{0.969} & \textbf{0.958}\\
			& & $\mathcal{N}(1,1)$ & \textbf{0.959} & 0.925 & 0.929 & 0.886 & \textbf{0.976} & \textbf{0.945} & 0.939 & 0.938 \\
			\hline\hline
	\end{tabular}}
	\caption{\textit{Model-based: Average coverage probabilities with the coverage level being $0.95$ with $t$-distributed system error $\omega_t$.}}
	\label{tab:n_dis}
\end{table}

Tables \ref{tab:t_dis} and \ref{tab:n_dis} show the average coverage probabilities with a coverage level of $0.95$ for different types of $W_n$ and different values of $m$ when the system error $\omega_t$ follows a normal or $t$ distribution. We can see that the coverage probability of the bootstrapping {PG} estimator approaches the nominal level of $95 \%$ as $m$ increases. The weight schemes and underlying dynamics jointly impact empirical coverage. Although no theoretical distributional assumption is made for random noise (only a finite second moment is required), dynamics with heavy-tailed errors ($t$ distribution) require a larger sample size of agents than normal cases to ensure empirical coverage. In addition, both Tables \ref{tab:t_dis} and \ref{tab:n_dis} show that the restriction to positive random variables could affect the efficiency: it tends to underestimate the actual confidence region.

\section{Conclusion} \label{sec:discussion}

In this paper, we propose a bootstrap procedure for statistical inference of the policy gradient for {LQ RL}. We establish the \textit{first} central limit theorem for the policy gradient estimator and the distributional consistency of the proposed bootstrap approach. Our experimental results demonstrated the effectiveness of the proposed method in ensuring coverage for noisy linear dynamic systems with {LQ RL} tasks. We note that for a broader class of {RL} "similar" to the {LQ RL} framework with stochastic dynamics, the proposed procedure could learn the inferential tasks. Although model-based {LQ RL} has been shown in the literature to outperform the other control models when the true system is linear, we often work with a finite sample size of agents and are uncertain whether the actual system is in a linear structure in the learning environment. We are interested in exploring the finite-sample properties and sensitivity of the proposed approach when the model is misspecified, e.g., there might be few or even single agents in the system dynamics, and considering extending the current work to more complex uncertain dynamics, which is more robust and scalable to model misspecification and small samples. Hypothesis testing is another aspect of statistical inference that is considered. Due to the distributional consistency of the proposed bootstrap technique, this can be developed under the same framework.

\begin{table}[htbp]\caption{Notation Summary} \footnotesize
	\begin{center}
		\begin{tabular}{r c p{11cm} } 
			\toprule
			$T$, $t$ & & Size of finite time horizon. $t$ indicates the $t$th time point \\
			$x_t \in \mathbb{R}^{d}$, $u_t \in \mathbb{R}^{p}$, $\omega_t \in \mathbb{R}^{d}$ & &  The state of the system, the action taken, the random noisy at time $t$. \\
			$\mathcal{P}_t$ && $x_t \sim \mathcal{P}$. \\
			$A \in \mathbb{R}^{d \times d}$, $B \in \mathbb{R}^{d \times p}$ & & The system (transition) matrices.\\
			$Q_t \in \mathbb{R}^{d \times d}$, $R_t \in \mathbb{R}^{p \times p}$ &  &  The parameterized matrices define the objectives at time $t$.\\
			$K_t \in \mathbb{R}^{d \times p}$,  $K_t^{\star} \in \mathbb{R}^{d \times p}$ && The policy at time $t$. The optimal policy at time $t$. \\
			$P_t^{\star}$, $P_t^{\boldsymbol{K}}$ && The optimal solution to the discrete algebraic Riccati equation at time $t$. The solution to the discrete algebraic Riccati equation at time $t$ for the policy gradient of $\boldsymbol{K}$. \\
			$E_t\left(K_t\right)$, $\Sigma_t$ && $E_t\left(K_t\right) :=\left(R_t+B^{\top} P_{t+1} B\right) K_t-B^{\top} P_{t+1} A$, $\Sigma_t := \mathbb{E}\left[x_t x_t^{\top}\right]$. \\
			$m$, $i$ && Size of the trajectory of the observation, $i$ indicates the $i$th trajectory.\\
			$\ell$, $n$, $\eta_\ell$ && $\ell$ indicates $\ell$th iteration of policy gradient descent. The maximum number of iterations. The stepsize of policy gradient at $\ell$th iteration. \\
			$K_t^{\ell}$, $\tilde{K}_t^{\ell}$ && The $\ell$th iterates of exact policy gradient descent at time $t$. The $\ell$th iterates of zeroth order policy gradient descent at time $t$. \\
			$U_t$ && The random matrix generated from a Uniform distribution with radius of $r$ at time $t$. \\
			$\bar{K}_t^n$, $\bar{\tilde{K}}_t^n$ &&  Polyak-Ruppert averaging iterate of exact policy gradient descent $\bar{K}_t^n=\frac{1}{n} \sum_{\ell=1}^n K_t^{\ell}$, and zeroth order policy gradient descent $\bar{\tilde{K}}_t^n=\frac{1}{n} \sum_{\ell=1}^n \tilde{K}_t^{\ell}$ \\
			$c(\cdot)$,  $\widehat{c}(\cdot)$ && The risk given  state and input $K$. The real observation of risk given an arbitrary state and input $K$. \\
			$W_\ell$ && The random variable for the Bootstrap procedure at $\ell$th iteration with $\mathbb{E}[W_\ell] = 1$ and $\operatorname{Var}[W_\ell] = 1$, $\forall \ell$. \\
			$\mathcal{B}$, $b$ && Size of Bootstrap iteration, $b$ indicates $b$th iteration. \\
			$\bar{K}_{t, (b)}^{n, \prime}$ &&  Polyak-Ruppert averaging iterate of exact policy gradient descent $\bar{K}_{t, (b)}^{n, \prime}=\frac{1}{n} \sum_{\ell=1}^n K_{t, (b)}^{n, \prime}$, for $b$th bootstrap iteration. \\
			$\bar{\tilde{K}}_{t, (b)}^{n, \prime}$  &&  Polyak-Ruppert averaging iterate of zeroth order policy gradient descent $\bar{\tilde{K}}_{t, (b)}^{n, \prime}=\frac{1}{n} \sum_{\ell=1}^n \tilde{K}_{t, (b)}^{n, \prime}$, for $b$th bootstrap iteration. \\
			$\boldsymbol{K}: = \{K_0, \cdots K_{T - 1}\}$ && The vector representation of the policy matrix.\\
			$\boldsymbol{K}^{\star}: = \{K_0^{\star}, \cdots K^{\star}_{T - 1}\}$ && The vector representation of the optimal policy matrix.\\
			$C(\boldsymbol{K})$ && The long-time risk (Objective function value) given policy $\boldsymbol{K}$. \\
			$\nabla_t C(\boldsymbol{K}^\ell)$,  $\widehat{\nabla}_t C(\boldsymbol{K}^\ell)$ && The exact gradient of $\boldsymbol{K}^\ell$ at $\ell$th iteration at time $t$. The zeroth order gradient of $\boldsymbol{K}^\ell$ at $\ell$th iteration at time $t$.  \\
			$\boldsymbol{U}_t$ && The vector representation of a random matrix generated from a Uniform distribution with radius of $r$ at time $t$. \\
			\bottomrule
		\end{tabular}
	\end{center}
	\label{tab:not}
\end{table}

\newpage

\bibliographystyle{plain}

\bibliography{reference.bib}

\begin{thebibliography}{10}

\bibitem{anastasiou2019normal}
Andreas Anastasiou, Krishnakumar Balasubramanian, and Murat~A Erdogdu.
\newblock Normal approximation for stochastic gradient descent via
  non-asymptotic rates of martingale clt.
\newblock In {\em Conference on Learning Theory}, pages 115--137. PMLR, 2019.

\bibitem{aastrom1995adaptive}
Karl~Johan {\AA}str{\"o}m.
\newblock Adaptive control.
\newblock In {\em Mathematical System Theory: The Influence of RE Kalman},
  pages 437--450. Springer, 1995.

\bibitem{basei2022logarithmic}
Matteo Basei, Xin Guo, Anran Hu, and Yufei Zhang.
\newblock Logarithmic regret for episodic continuous-time linear-quadratic
  reinforcement learning over a finite-time horizon.
\newblock {\em Journal of Machine Learning Research}, 23(178):1--34, 2022.

\bibitem{becker1985adaptive}
Arthur Becker, P~Kumar, and Ching-Zong Wei.
\newblock Adaptive control with the stochastic approximation algorithm:
  Geometry and convergence.
\newblock {\em IEEE Transactions on Automatic Control}, 30(4):330--338, 1985.

\bibitem{bertsekas2012dynamic}
Dimitri Bertsekas.
\newblock {\em Dynamic programming and optimal control: Volume I}, volume~4.
\newblock Athena scientific, 2012.

\bibitem{bradtke1992reinforcement}
Steven Bradtke.
\newblock Reinforcement learning applied to linear quadratic regulation.
\newblock {\em Advances in Neural Information Processing Systems}, 5, 1992.

\bibitem{carmona2019linear}
Ren{\'e} Carmona, Mathieu Lauri{\`e}re, and Zongjun Tan.
\newblock Linear-quadratic mean-field reinforcement learning: convergence of
  policy gradient methods.
\newblock {\em arXiv preprint arXiv:1910.04295}, 2019.

\bibitem{chen2021statistical}
Haoyu Chen, Wenbin Lu, and Rui Song.
\newblock Statistical inference for online decision making via stochastic
  gradient descent.
\newblock {\em Journal of the American Statistical Association},
  116(534):708--719, 2021.

\bibitem{chen2020statistical}
Xi~Chen, Jason~D Lee, Xin~T Tong, and Yichen Zhang.
\newblock Statistical inference for model parameters in stochastic gradient
  descent.
\newblock 2020.

\bibitem{chen2020robust}
Xi~Chen and Wen-Xin Zhou.
\newblock Robust inference via multiplier bootstrap.
\newblock 2020.

\bibitem{dean2020sample}
Sarah Dean, Horia Mania, Nikolai Matni, Benjamin Recht, and Stephen Tu.
\newblock On the sample complexity of the linear quadratic regulator.
\newblock {\em Foundations of Computational Mathematics}, 20(4):633--679, 2020.

\bibitem{ding2020challenges}
Zihan Ding and Hao Dong.
\newblock Challenges of reinforcement learning.
\newblock {\em Deep Reinforcement Learning: Fundamentals, Research and
  Applications}, pages 249--272, 2020.

\bibitem{dulac2021challenges}
Gabriel Dulac-Arnold, Nir Levine, Daniel~J Mankowitz, Jerry Li, Cosmin
  Paduraru, Sven Gowal, and Todd Hester.
\newblock Challenges of real-world reinforcement learning: definitions,
  benchmarks and analysis.
\newblock {\em Machine Learning}, 110(9):2419--2468, 2021.

\bibitem{fallah2021convergence}
Alireza Fallah, Kristian Georgiev, Aryan Mokhtari, and Asuman Ozdaglar.
\newblock On the convergence theory of debiased model-agnostic
  meta-reinforcement learning.
\newblock {\em Advances in Neural Information Processing Systems},
  34:3096--3107, 2021.

\bibitem{fang2018online}
Yixin Fang, Jinfeng Xu, and Lei Yang.
\newblock Online bootstrap confidence intervals for the stochastic gradient
  descent estimator.
\newblock {\em Journal of Machine Learning Research}, 19(78):1--21, 2018.

\bibitem{fazel2018global}
Maryam Fazel, Rong Ge, Sham Kakade, and Mehran Mesbahi.
\newblock Global convergence of policy gradient methods for the linear
  quadratic regulator.
\newblock In {\em International conference on machine learning}, pages
  1467--1476. PMLR, 2018.

\bibitem{giegrich2024convergence}
Michael Giegrich, Christoph Reisinger, and Yufei Zhang.
\newblock Convergence of policy gradient methods for finite-horizon exploratory
  linear-quadratic control problems.
\newblock {\em SIAM Journal on Control and Optimization}, 62(2):1060--1092,
  2024.

\bibitem{guo2023reinforcement}
Xin Guo, Anran Hu, and Yufei Zhang.
\newblock Reinforcement learning for linear-convex models with jumps via
  stability analysis of feedback controls.
\newblock {\em SIAM Journal on Control and Optimization}, 61(2):755--787, 2023.

\bibitem{hambly2021policy}
Ben Hambly, Renyuan Xu, and Huining Yang.
\newblock Policy gradient methods for the noisy linear quadratic regulator over
  a finite horizon.
\newblock {\em SIAM Journal on Control and Optimization}, 59(5):3359--3391,
  2021.

\bibitem{hao2019bootstrapping}
Botao Hao, Yasin Abbasi~Yadkori, Zheng Wen, and Guang Cheng.
\newblock Bootstrapping upper confidence bound.
\newblock {\em Advances in Neural Information Processing Systems}, 32, 2019.

\bibitem{huang2020model}
Qingyan Huang.
\newblock Model-based or model-free, a review of approaches in reinforcement
  learning.
\newblock In {\em 2020 International Conference on Computing and Data Science
  (CDS)}, pages 219--221. IEEE, 2020.

\bibitem{ji2019improved}
Kaiyi Ji, Zhe Wang, Yi~Zhou, and Yingbin Liang.
\newblock Improved zeroth-order variance reduced algorithms and analysis for
  nonconvex optimization.
\newblock In {\em International conference on machine learning}, pages
  3100--3109. PMLR, 2019.

\bibitem{keskin2018incomplete}
N~Bora Keskin and Assaf Zeevi.
\newblock On incomplete learning and certainty-equivalence control.
\newblock {\em Operations Research}, 66(4):1136--1167, 2018.

\bibitem{lai1982iterated}
Tze~Leung Lai and Herbert Robbins.
\newblock Iterated least squares in multiperiod control.
\newblock {\em Advances in Applied Mathematics}, 3(1):50--73, 1982.

\bibitem{lee2024fast}
Sokbae Lee, Yuan Liao, Myung~Hwan Seo, and Youngki Shin.
\newblock Fast inference for quantile regression with tens of millions of
  observations.
\newblock {\em Journal of Econometrics}, page 105673, 2024.

\bibitem{li2021statistical}
Xiang Li, Jiadong Liang, Xiangyu Chang, and Zhihua Zhang.
\newblock Statistical estimation and inference via local sgd in federated
  learning.
\newblock {\em arXiv preprint arXiv:2109.01326}, 2021.

\bibitem{li2021unifying}
Ye~Li, Hong Xie, Yishi Lin, and John~CS Lui.
\newblock Unifying offline causal inference and online bandit learning for data
  driven decision.
\newblock In {\em Proceedings of the Web Conference 2021}, pages 2291--2303,
  2021.

\bibitem{littman1996generalized}
Michael~L Littman and Csaba Szepesv{\'a}ri.
\newblock A generalized reinforcement-learning model: Convergence and
  applications.
\newblock In {\em ICML}, volume~96, pages 310--318, 1996.

\bibitem{lockwood2022review}
Owen Lockwood and Mei Si.
\newblock A review of uncertainty for deep reinforcement learning.
\newblock In {\em Proceedings of the AAAI Conference on Artificial Intelligence
  and Interactive Digital Entertainment}, volume~18, pages 155--162, 2022.

\bibitem{lovatto2021gradient}
{\^A}ngelo~Greg{\'o}rio Lovatto, Thiago~Pereira Bueno, and Leliane~Nunes
  de~Barros.
\newblock Gradient estimation in model-based reinforcement learning: a study on
  linear quadratic environments.
\newblock In {\em Brazilian Conference on Intelligent Systems}, pages 33--47.
  Springer, 2021.

\bibitem{mania2019certainty}
Horia Mania, Stephen Tu, and Benjamin Recht.
\newblock Certainty equivalence is efficient for linear quadratic control.
\newblock {\em Advances in Neural Information Processing Systems}, 32, 2019.

\bibitem{meadows2022linear}
Olusesi~Ayobami Meadows, Arthur Rodriguez, and Ahmed~Tijani Salawudeen.
\newblock A linear quadratic regulator based speed control for
  remote-controlled racing cars.
\newblock In {\em 2022 IEEE Nigeria 4th International Conference on Disruptive
  Technologies for Sustainable Development (NIGERCON)}, pages 1--5. IEEE, 2022.

\bibitem{nikolakakis2022black}
Konstantinos Nikolakakis, Farzin Haddadpour, Dionysis Kalogerias, and Amin
  Karbasi.
\newblock Black-box generalization: Stability of zeroth-order learning.
\newblock {\em Advances in neural information processing systems},
  35:31525--31541, 2022.

\bibitem{polyak1992acceleration}
Boris~T Polyak and Anatoli~B Juditsky.
\newblock Acceleration of stochastic approximation by averaging.
\newblock {\em SIAM Journal on Control and Optimization}, 30(4):838--855, 1992.

\bibitem{pong2018temporal}
Vitchyr Pong, Shixiang Gu, Murtaza Dalal, and Sergey Levine.
\newblock Temporal difference models: Model-free deep rl for model-based
  control.
\newblock {\em arXiv preprint arXiv:1802.09081}, 2018.

\bibitem{ramprasad2023online}
Pratik Ramprasad, Yuantong Li, Zhuoran Yang, Zhaoran Wang, Will~Wei Sun, and
  Guang Cheng.
\newblock Online bootstrap inference for policy evaluation in reinforcement
  learning.
\newblock {\em Journal of the American Statistical Association},
  118(544):2901--2914, 2023.

\bibitem{rathnam2023unintended}
Sarah Rathnam, Sonali Parbhoo, Weiwei Pan, Susan Murphy, and Finale
  Doshi-Velez.
\newblock The unintended consequences of discount regularization: Improving
  regularization in certainty equivalence reinforcement learning.
\newblock In {\em International Conference on Machine Learning}, pages
  28746--28767. PMLR, 2023.

\bibitem{razin2024implicit}
Noam Razin, Yotam Alexander, Edo Cohen-Karlik, Raja Giryes, Amir Globerson, and
  Nadav Cohen.
\newblock Implicit bias of policy gradient in linear quadratic control:
  Extrapolation to unseen initial states.
\newblock {\em arXiv preprint arXiv:2402.07875}, 2024.

\bibitem{ruppert1988efficient}
David Ruppert.
\newblock Efficient estimations from a slowly convergent robbins-monro process.
\newblock Technical report, Cornell University Operations Research and
  Industrial Engineering, 1988.

\bibitem{samsonov2024gaussian}
Sergey Samsonov, Eric Moulines, Qi-Man Shao, Zhuo-Song Zhang, and Alexey
  Naumov.
\newblock Gaussian approximation and multiplier bootstrap for polyak-ruppert
  averaged linear stochastic approximation with applications to td learning.
\newblock In {\em The Thirty-eighth Annual Conference on Neural Information
  Processing Systems}, 2024.

\bibitem{shi2023dynamic}
Chengchun Shi, Xiaoyu Wang, Shikai Luo, Hongtu Zhu, Jieping Ye, and Rui Song.
\newblock Dynamic causal effects evaluation in a/b testing with a reinforcement
  learning framework.
\newblock {\em Journal of the American Statistical Association},
  118(543):2059--2071, 2023.

\bibitem{sun2021learning}
Yue Sun and Maryam Fazel.
\newblock Learning optimal controllers by policy gradient: Global optimality
  via convex parameterization.
\newblock In {\em 2021 60th IEEE Conference on Decision and Control (CDC)},
  pages 4576--4581. IEEE, 2021.

\bibitem{szpruch2024optimal}
Lukasz Szpruch, Tanut Treetanthiploet, and Yufei Zhang.
\newblock Optimal scheduling of entropy regularizer for continuous-time
  linear-quadratic reinforcement learning.
\newblock {\em SIAM Journal on Control and Optimization}, 62(1):135--166, 2024.

\bibitem{umenberger2019robust}
Jack Umenberger, Mina Ferizbegovic, Thomas~B Sch{\"o}n, and H{\aa}kan
  Hjalmarsson.
\newblock Robust exploration in linear quadratic reinforcement learning.
\newblock {\em Advances in Neural Information Processing Systems}, 32, 2019.

\bibitem{vlahakis2019distributed}
Eleftherios~E Vlahakis, Leonidas~D Dritsas, and George~D Halikias.
\newblock Distributed lqr design for identical dynamically coupled systems:
  Application to load frequency control of multi-area power grid.
\newblock In {\em 2019 IEEE 58th Conference on Decision and Control (CDC)},
  pages 4471--4476. IEEE, 2019.

\bibitem{wagenmaker2022first}
Andrew~J Wagenmaker, Yifang Chen, Max Simchowitz, Simon Du, and Kevin Jamieson.
\newblock First-order regret in reinforcement learning with linear function
  approximation: A robust estimation approach.
\newblock In {\em International Conference on Machine Learning}, pages
  22384--22429. PMLR, 2022.

\bibitem{wang2020residual}
Chi-Hua Wang, Yang Yu, Botao Hao, and Guang Cheng.
\newblock Residual bootstrap exploration for bandit algorithms.
\newblock {\em arXiv preprint arXiv:2002.08436}, 2020.

\bibitem{wang2021exact}
Feicheng Wang and Lucas Janson.
\newblock Exact asymptotics for linear quadratic adaptive control.
\newblock {\em Journal of Machine Learning Research}, 22(265):1--112, 2021.

\bibitem{wang2020continuous}
Haoran Wang and Xun~Yu Zhou.
\newblock Continuous-time mean--variance portfolio selection: A reinforcement
  learning framework.
\newblock {\em Mathematical Finance}, 30(4):1273--1308, 2020.

\bibitem{wang2018stochastic}
Yining Wang, Simon Du, Sivaraman Balakrishnan, and Aarti Singh.
\newblock Stochastic zeroth-order optimization in high dimensions.
\newblock In {\em International conference on artificial intelligence and
  statistics}, pages 1356--1365. PMLR, 2018.

\bibitem{yang2020leveraging}
Haojun Yang, Kuan Zhang, Kan Zheng, and Yi~Qian.
\newblock Leveraging linear quadratic regulator cost and energy consumption for
  ultrareliable and low-latency iot control systems.
\newblock {\em IEEE Internet of Things Journal}, 7(9):8356--8371, 2020.

\bibitem{yang2019provably}
Zhuoran Yang, Yongxin Chen, Mingyi Hong, and Zhaoran Wang.
\newblock Provably global convergence of actor-critic: A case for linear
  quadratic regulator with ergodic cost.
\newblock {\em Advances in Neural Information Processing Systems}, 32, 2019.

\bibitem{zhu2023online}
Wanrong Zhu, Xi~Chen, and Wei~Biao Wu.
\newblock Online covariance matrix estimation in stochastic gradient descent.
\newblock {\em Journal of the American Statistical Association},
  118(541):393--404, 2023.

\bibitem{zhu2024uncertainty}
Yi~Zhu, Jing Dong, and Henry Lam.
\newblock Uncertainty quantification and exploration for reinforcement
  learning.
\newblock {\em Operations Research}, 72(4):1689--1709, 2024.

\end{thebibliography}

\end{document}